\documentclass[12pt]{article}
\usepackage{amssymb}
\usepackage{pb-diagram}

\title{Coactions of Hopf-$C^*$-algebras and equivariant $E$-theory }
\author{Radu POPESCU \thanks{ Supported by the FCT grant no. SFRH/BPD/20306/2004,
at the Instituto Superior Tecnico in Lisbon.}}
\date{}

\def\a{\delta_A}
\def\ab{\delta_{A\o \K}}
\def\b{\delta_B}
\def\bb{\delta_{B\o \K}}
\def\c{\delta_C}

\def\ss{\bar S\mathrm{-alg}}
\def\U{\mathfrak{A}}

\def\C{$C^*$-algebra}
\def\M{\tilde{M}}
\def\d{\delta}
\def\>{\rightarrow}

\def\f{\varphi}
\def\L{{\cal L}}
\def\K{{\cal K}}

\def\o{\otimes}
\def\r{\rtimes_r}
\def\im{\mathrm{im}}
\def\id{\mathrm{id}}

\newtheorem{teo}{Theorem}[section]
\newtheorem{lema}[teo]{Lemma}
\newtheorem{defi}[teo]{Definition}
\newtheorem{prop}[teo]{Proposition}
\newtheorem{coro}[teo]{Corollary}
\newtheorem{rmq}[teo]{Remark}

\newcommand{\ma}[4]{\left(\begin{array}{cc}
#1 & #2\\
#3 & #4
\end{array}\right)}

\begin{document}

\maketitle

\begin{abstract}
We define and study an equivariant $E$-theory with respect to coactions of Hopf
$C^*$-algebras;  we prove the Baaj-Skandalis duality in this setting. We show that the
corresponding equivariant $KK$-theory of Baaj and Skandalis enjoys an universal property.
In the appendix, we look at the different ways of expressing equivariant stability for a
functor, and prove an equivariant Brown-Green-Rieffel stabilization result.
\end{abstract}


\tableofcontents

\section{Introduction}

Assume that $A$ is a $C^*$-algebra acted upon by a locally compact abelian group $G$; the
cross product $A\rtimes G$ is endowed with an action of the dual group of $G$, $\hat{G}$,
and the double cross product $A\rtimes G \rtimes \hat G$ is isomorphic to $A\o
\K(L^2(G))$. This isomorphism, the Takesaki-Takay duality, was generalized to
noncommutative groups, and roughly states that the cross product $A\rtimes_r G$ is
endowed with a coaction of $G$, and a similar duality holds by taking the cross product
by the coaction of $G$. When $G$ is abelian, coactions of $G$ correspond to actions of
the dual group $\hat G$, so one recovers the same isomorphism.

Hopf $C^*$-algebras provide a common framework for actions and coactions of groups and
their cross products, and for duality results. It also allows us  to prove duality beyond
the case of groups, in the realm of quantum groups.
 In the
general setting, Hopf $C^*$-algebras were considered in $K$-theory by Baaj and Skandalis
who defined and studied an equivariant $KK$-theory with respect to coactions of Hopf
$C^*$-algebras.  This theory  is a generalization of Kasparov's equivariant $KK$-theory
for groups: when the Hopf $C^*$-algebra in question is $C_0(G)$, the Hopf algebra
associated to a group $G$, there is an identification between $KK^{C_0(G)}$ and $KK^G$.

One of the main properties of this theory is the Baaj-Skandalis  isomorphism  which
identifies  the group $KK^G(A,B)$ with the group  $KK^{C^*_r(G)}(A\rtimes_r G, B
\rtimes_r G)$; there is also a dual isomorphism.

 This Baaj-Skandalis isomorphism is a valuable tool in the study of the
$K$-theory of cross product algebras, and recently J. Cuntz used it in connection to the
Baum-Connes conjecture which predicts that the $K$-theory of the cross product
$A\rtimes_r G$ is isomorphic, by means of an index map, to a group of generalized
elliptic operators with coefficients in $A$. In \cite{Cu|BC} Cuntz reinterpreted the
above conjecture for discrete groups as follows: using the Baaj-Skandalis isomorphism he
defined a finitely supported $K$-theory, $K^{fin}_{*}(A\rtimes_r G)$, for the cross
product algebras $A\rtimes_r G$, and reformulated the Baum-Connes assembly map as an
application
$$\mu: K^{fin}_*(A\rtimes_r G)\>K_*(A\rtimes_r G)$$
between finitely supported and ordinary $K$-theory of the cross product algebras
$A\rtimes_r G$.

In this article we define the equivariant $E$-theory for coactions of Hopf
$C^*$-algebras. Unless stated otherwise, we use the minimal tensor product of
$C^*$-algebras and the reduced cross product.

$E$-theory, which is due to Connes and Higson, associates to a pair of $C^*$-algebras
$(A,B)$, an abelian group $E(A,B)$. It is a contravariant functor in the first variable
and a covariant one in the second, and there is a product $E(A,B)\times E(B,C)\> E(A,C)$
which allows us to see it  as a category.  All these properties are also shared by the
Kasparov's $KK$-theory of which $E$-theory is a variant : there is  a natural
transformation from $KK$ to $E$-theory, and various results explain the relations and the
differences between these theories.

Equivariant $E$-theory with respect to group actions was defined by Guentner, Higson and
Trout (\cite{GHT}) and was used by Higson and Kasparov in their proof of the Baum-Connes
conjecture for amenable groups (\cite{HK2}).

An important feature of $E$-theory is that, contrary to $KK$-theory, it does not behave
well with respect to minimal tensor product of $C^*$-algebras or reduced cross products
by groups, mainly because these operations are not exact, as explain in Section 6.3 of
\cite{GH|LNM}, based on the argument from \cite{S}. But $E$-theory is based on the
following simple to state notion, introduced by Connes and Higson in \cite{CH}:

\begin{defi} \em Let $T=[1,+\infty[$. An asymptotic morphism between two $C^*$-algebras $A$
and $B$ is a family of maps $(\f_t)_{t\in T} :A\> B$ such that
\begin{itemize}
\item[(i)]
$\f_t(a+\lambda b)-\f_t(a)-\lambda\f_t(b)\>0$
\item[(ii)]
$\f_t(ab)-\f_t(a)\f_t(b)\>0$
\item[(iii)]
$\f_t(a)^*-\f_t(a^*)\>0$
\end{itemize}
 as $t\>+ \infty$ and such that for every $a\in A$ , $t\mapsto \f_t(a)$ is continuous.
\end{defi}

Two asymptotic morphisms $(\f_t)_{t\in T}$ and $(\f'_t)_{t\in T}$  are (asymptotically)
equivalent if for every $a\in A$, $(\f_t-\f'_t)(a)\>0$ as $t\>+\infty$. Two asymptotic
morphisms $\f_t^0,\f_t^1:A\>B$ are {\it homotopic} if they are obtained from an
asymptotic morphism $\Phi:A\>C([0,1],B)$ by evaluation in $0$ and in $1$, respectively.
Homotopy is an equivalence relation, denote by $\{A,B\}$ the homotopy classes of
asymptotic morphisms between $A$ and $B$.

A different picture of this is obtained as follows.

\begin{defi} \em
Let $A$ be a $C^*$-algebra, the  asymptotic algebra of $A$, $\U A =C_b(T,A)/C_0(T,A)$ is
the quotient of the $C^*$-algebra of bounded continuous functions on $T$ with values in
$A$ by the ideal of functions vanishing at infinity. Denote by $\alpha_A : A\> \U A$ the
map that associates to an element $a\in A$ the class of the constant map $t\mapsto a$.
\end{defi}

There is a one-to-one correspondence between equivalence classes of asymptotic morphisms
and *-homomorphisms $\f : A\> \U B$. We shall follow \cite{GHT} and use this point of
view, rather then the initial approach from \cite{CH}. Let $\f:A\>B$ be a
$*$-homomorphism between the $C^*$-algebras $A$ and $B$; composition with $\f$ induces a
map $\U \f:\U A\> \U B$, so one can regard $\U$ as a functor on the category of
$C^*$-algebras. We use the notation $\U^n$ for the $n$-fold composition of $\U$ with
itself for $n \geq 1$ and the identity functor for $n=0$.

A summary of this article is as follows.

We start the first section by recalling some facts and by fixing some notations about
multiplier algebras, Hopf $C^*$-algebras $S$, and coactions of them on $C^*$-algebras.
Details can be found in \cite{BS} and in \cite{EKQR}. We define the notion of an
$S$-asymptotic morphism and $\U_S A$, the subalgebra of $S$-continuous elements of $\U
A$. We show that the asymptotic morphism associated by the Connes-Higson construction to
an $S$-equivariant extension is equivariant in the sense that we had defined.

 In the second part
we construct a category which combines homotopy classes of equivariant asymptotic
morphism  with exterior equivalence of coactions. We keep track of this using cocycles
for the coactions involved.

Then we define the equivariant $E$-theory with respect to coactions of Hopf
$C^*$-algebras and prove its main properties:existence of a product, six-exact sequences,
tensor product. We also show that, when the Hopf algebra is associated to the action of a
locally compact group, i.e., when $S=C_0(G)$, we can identify our $E^{C_0(G)}$ with the
equivariant $E$-theory for groups, $E^G$, previously defined by Gueunter, Higson and
Trout in \cite{GHT}. We obtain in this way an alternative description of their theory. We
prove the Baaj-Skandalis isomorphism in $E$-theory.

 In the fourth section we study the equivariant $KK$-theory from the point of view of
its universal property and derive the existence of a natural transformation from $KK^S$
to $E^S$.

In the appendix we prove an equivariant version of the a theorem of Brown, Green,
Rieffel, which states that two Morita equivalent separable $C^*$-algebras are stable
isomorphic. As a consequence one obtains that different ways of expressing stability of a
functor are equivalent.

This article was written while the author was a postdoc in Muenster, Cardiff, and Lisbon.

\section{Equivariant asymptotic morphisms}

Let $A$ be a $C^*$-algebra, denote by $M(A)$  the associated multiplier algebra, the
biggest $C^*$-algebra which contains $A$ as an essential ideal. If $J$ is an indeal in
$A$, it follows that there is a restriction map from $M(J)$ to $M(A)$; this map is
neither injective nor surjective, in general.

$M(A)$ is the completion of $A$ with respect to the strict topology, the topology
generated by the seminorms $p_c(a) = \|ac\| + \|ca\|$, indexed by elements $c$ in $A$.
Recall that a morphism $\f:A\>B$ is {\it nondegenerate} if there is an approximate unit
$(e_i)_i$ for $A$ such that $\f(e_i)$ strictly converges to $1_{M(B)}$. Such a morphism
$\f$ extends to a map between the multiplier algebras, a map that we denote again by
$\f:M(A)\>M(B)$.

Let $A$ and $S$ be \C s, denote by
$$ \M(A\o S)=\{ m\in M(A\o S): m(1\o S)+ (1\o S)m \subset A\o S\}, $$
the $S$-multipliers of $A\o S$; this is a closed *-subalgebra of $M(A\o S)$.

 The $S$-strict topology of $\M(A\o S)$ is the
locally convex topology generated by the seminorms $p_s(m)= \|m(1\o s)\| + \|(1\o s)m\|$
indexed by  $s\in S$.

The relevance of this topology for $\M(A\o S)$ was observed in \cite{EKQR}, Proposition
A.4, where it is proved that:
\begin{lema} $\M(A\o S)$ is the $S$-strict
completion of $A\o S$.
\end{lema}
{\bf Examples} 1. If $S$ is unital, then $\M(A\o S)= A\o S$.\\
2. Assume that $A$ is unital, then $\M(A\o S)=M(A\o S)$.\\
3. When $S=C_0(X)$ the continuous functions vanishing at infinity on a locally compact
Hausdorff space $X$, then $\M(A\o C_0(X))=C_b(X,A)$, the algebra of bounded continuous
functions with values in $A$. Note that the multiplier algebra $M(A\o C_0(X))$ is
$C_b(X,M^s(A)) $ the algebra
of bounded strictly continuous functions, a much larger algebra in general.\\

Consider \C s, $A,B,C,D$ and *-homomorphisms $\f:A\>M(C)$ and $\psi:B\> M(D)$. Suppose
that $\psi$ is nondegenerate, then $\f\o\psi : A\o B\> M(C\o D) $ extends to $\M(A\o B)$
and when $\f(A)\subset C$, the image of this morphism sits inside $\M(C\o D)$ and the
restricted morphism, denoted again $\f\o\psi : \M(A\o B)\> \M(C\o D)$ is continuous from
the $B$-strict topology of $\M(A\o B)$ to the $D$-strict topology of $\M(C\o D)$.

In particular  any *-homomorphism  $\f:A\>M(B)$ gives rise to a map $\f\o \id_S :\M(A\o
S)\>M(B\o S)$.

In general, we use the same symbol for a $*$-homomorphism and for its extension to the
$S$-multiplier algebras.

\begin{defi} \em A Hopf \C \ is a \C \ $S$ endowed with a
nondegenerate *-homomorphism $\d_S : S\> \M(S\o S)$ which verify the following
comultiplication condition $$(\id_S\o \d_S)\circ \d_S = (\d_S\o \id_S)\circ \d_S.$$ A
coaction of a Hopf \C \ $S$ on a \C \ $A$ is a nondegenerate *-homomorphism $\a : A\>
\M(A\o S)$ such that
$$(\id_A\o \d_S)\circ \a = (\a\o \id_S)\circ \a .$$
\end{defi}

This condition is equivalent to the commutativity of the following diagram:

$$
\begin{diagram}
\node{A}\arrow{e,t}{\d_A}\arrow{s,l}{\d_A}\node{\M(A\o S)}\arrow{s,r}{\id_A\o \d_S}\\
\node{\M(A\o S)}\arrow{e,t}{\d_A \o \id_S} \node{M(A \o S \o S);}
\end{diagram}
$$
we denote by $\a^2:=(\a\o \id_S) \a :A\>M(A\o S\o S)$ the above morphism.

In this case we shall say that $A$ is an $S$-algebra. An element $a\in A$ is {\it fixed}
by the coaction $\a$ if $\d_A(a)=a\o 1$.

 Let $A$ and $B$ be $S$-algebras, a *-homomorphism  $\f:A\>M(B)$ is $S$-equivariant
 (or an $S$-morphism) if $(\f\o\id_S)\d_A=\d_B \f$.

\begin{defi} \em A morphism between two  Hopf \C s $S$ and $S'$ is a
nondegenerate morphism $\f: S\> M(S')$ such that the following diagram commutes
$$
\begin{diagram}
\node{S}\arrow{e,t}{\f}\arrow{s,l}{\d_S}\node{M(S')}\arrow{s,r}{\d_{S'}} \\
\node{\M(S\o S)}\arrow{e,t}{\f \o \f} \node{M(S' \o S').}
\end{diagram}
$$
\end{defi}

 Such a morphism $\f$ allows us to define an coaction of $S'$ on any $S$-algebra $A$ by
$\d'_A(a)=(\id_A\o \f) \d_A$. \\
\\{\bf Examples} 1. Any \C \ $A$ has a {\it trivial} Hopf structure with
comultiplication
defined by $\delta(a)=a\o 1$.\\
2.{\it Actions of groups.}  Let $G$ be a locally compact group (or just a semigroup).
Define $\d:C_0(G)\>C_b(G,C_0(G))=\M(C_0(G)\o C_0(G))$  by $\d(f)(s,t)=f(st)$.
Associativity of the product of $G$ translates into the comultiplication property of $\d$
and it is easy to see that $C_0(G)$ is a Hopf $C^*$-algebra. Strongly continuous actions
of $G$ on an algebra $A$ correspond to injective coactions of $C_0(G)$ on $A$ in the
following way: let $\alpha$ be an action of $G$  on $A$, then the corresponding coaction
$\d_A:A\>C_b(G,A)$ is given by $\d_A(a)=g \mapsto\alpha_g(a)$
for all $a\in A$. \\
3.{\it Coactions of groups.} Let $G$ be a locally compact group. The reduced
$C^*$-algebra  of $G$, $C^*_r(G)$, is a Hopf \C \ with comultiplication $\delta_G$ given
by $\d (x)= \int x(g)(\lambda_g \o \lambda_g)dg$ for an element $x\in C^*_r(G), x=\int
x(g)\lambda_g dg$.

Define the unitary $W$ on $L^2(G)\o L^2(G)=L^2(G\times G)$ by $Wf(s,t)=f(s,s^{-1}t)$,
then the above coproduct can be written as $\d(x)=W(x\o 1)W^*$ from which follows that
$\d$ is nondegenerate.

A coaction of this Hopf $C^*$-algebra is called a coaction of $G$. If $G$ is commutative,
then the Fourier transform provides an isomorphism of Hopf \C s between $C^*_r(G)$ and
$C_0(\hat G)$. A coaction of $G$ amounts to an action of the dual group.

Let $A$ be a $G$-algebra,  the cross product $A\rtimes_r G$ is endowed with a coaction of
$G$, $\d_{A\rtimes_r G}:A\rtimes_r G\> \M(A\rtimes_r G \o C^*_r(G))$ given by
$\d_{A\rtimes_r G}(\int a(s)u_s ds)= \int (a(s)\o 1)(u_s\o \lambda _s) ds$.

Suppose that $A$ is endowed with a coaction of a {\it discrete} group $G$. For every
$g\in G$ its spectral subspace  is $A_g=\{a\in A \ \ \mbox{such that}\ \
\d_A(a)=a\o\lambda_g\}$, and $A$ is the closed linear span of the family $(A_g)_{g\in
G}$. The coaction condition translates into the fact that $(A_g)_{g\in G}$ is a
$G$-grading and to any $G$-grading there is an associated coaction.\\
4.{\it Multiplicative Unitaries.} The operator $W$ is an example of a multiplicative
unitary, introduced by Baaj and Skandalis in \cite{BS|unitmult}, i.e., an unitary $W\in
L(H\o H)$ which fulfills the pentagonal relation $W_{12}W_{13}W_{23}=W_{23}W_{12}$. Under
supplementary assumptions, associated to it there are two dual Hopf \C s and cross
product algebras constructions, further extending the results from actions and coactions
of groups.

\begin{lema}If the following sequence is exact
$$0\>J\o S\stackrel{i}{\>}B\o S \stackrel{p}{\>}A\o S\>0$$
then
$$ 0\>\M(J\o S)\stackrel{\bar{i}}{\>}\M(B\o S)\stackrel{\bar{p}}{\>}
\M(A\o S)\>0$$ is also exact.
\end{lema}
{\bf Proof.} For the injectivity, $i$ is the inclusion of $J\o S$ as an ideal inside $B\o
S$, $\bar{i}$ is then the inclusion of the corresponding $S$-completions. Equivalently,
$\im\ \bar{i}=\overline{\im\ i}^{S}$, the $S$-completion of the image of $i$.

We prove now that $\ker\bar{p}=\overline{\ker p}^S$. Take $b\in\ker\bar{p}$ and $s_i$ an
approximate unit for $S$; the sequence $b_n=b(1\o s_n)$ verifies $b_n\in B\o S \ ,
b_n\in\ker p \textrm{ as }\bar{p}(b(1\o s_n))=\bar{p}(b)(1\o s_n)=0$, and also $b_n\>b$
strictly, hence $b\in\overline{\ker p}^S$. Thus $\ker\bar{p}\subset\overline{\ker p}^S$.
The other inclusion follows from the strict continuity of $\bar{p}$. The initial short
exact sequence is exact in the middle thus
$$\im\ \bar{i}=\overline{\im  i}^S=\overline{\ker p}^S=\ker \bar{p}.$$

For the surjectivity of $\bar{p}:\M(B\o S)\> \M(A\o S)$ one can adapt the proof of the
noncommutative Tietze extension theorem \cite[Proposition 6.8]{L}, which states that the
corresponding map $\tilde{p}:M(B\o S)\> M(A\o S)$ is surjective. The proof quoted uses
the fact the multiplier algebras are strict completions of their corresponding algebras.
This proofs adapts to completions in the $S$-strict topology, that is, to $S$-multiplier
algebras. We outline it here.

First observe that a norm bounded sequence $(x_n)_n$ of elements of $A\o S$ converges
into the $S$-strict topology if the sequences $(x_n(1\o h))_n$ and $((1\o h)x_n)_n$
converge in the norm of $A\o S$, for some strictly positive element $h\in S$, cf.
\cite[Lemma 2.3.6.]{Wegge}.

Take $x\in \M(A\o S)$; we are looking for a lift $y\in \M(B\o S)$. If suffices to take
elements $x$ with $0\leq x \leq 1$; put $x_n:= x^{1/2}(1\o s_n)x^{1/2}$, where $(s_n)_n$
denotes an approximate unit for $S$. One has $(x_n)_n\in A\o S$ and $(x_n)_n$ is
increasing and converges to $x$ in the $S$-strict topology.

Choose $h$ a strictly positive element of $S$; then $((1\o h)x_n(1\o h))_n$ converges in
norm to $(1\o h)x(1\o h)$, thus, by eventually passing to a subsequence , we can assume
that
$$\|(x_{n+1}-x_n)(1\o h)\|< 2^{-n} \ \mbox{for all}\ n \ \mbox{and with}\ x_1=x_2=0.$$
By a lifting argument, there is a sequence $y_n\in B\o S$ such that
$$\|y_n\|\leq 1, \ \bar p(y_n)=x_n, \  \mbox{and }\ \|(y_{n+1}-y_n)(1\o h)\|< 2^{-n}.$$
The sequence $(y_n)_n$ converges in the $S$-strict topology of $\M(B \o S)$ to an element
$y$ such that $\bar p(y)=x$.
$\Box$\\
\\
Exactness plays a crucial role in the following.  Recall that a $C^*$-algebra $D$ is {\it
exact} if for every short exact sequence of $C^*$-algebras $$ 0\>J\>B\>A\>0$$ the
 sequence associated by taking minimal tensor product with $D$, i.e., if
$$ 0\>J\o S\>B\o S\>A\o S\>0$$
is also exact.

Related, a group $G$ is {\it exact}, if the cross product functor $\cdot \rtimes_r G$, is
exact, which means that for every short exact sequence of $G$-algebras $ 0\>J\>B\>A\>0$,
the sequence $ 0\>J\rtimes_r G\>B\rtimes_r G\>A\rtimes_r G\>0$ is also exact. In
particular, this implies that the reduced \C \ of  the group $G$ is exact; in the
discrete case this condition is also sufficient (\cite{KirWass|Exact groups}).

Note that if we replace the minimal tensor product with the maximal one, this condition
is always satisfied.

We shall consider only coactions of Hopf $C^*$-algebras $S$ which are exact as
$C^*$-algebras. This framework includes the cases when $S$ is nuclear, like in the case
of action by a group, when $S$ is commutative. It also include the case of coactions by
an exact group $G$.

We have chosen to work into the realm of minimal tensor products, but there is also a
theory using the maximal one, i.e., when working with coactions given by morphisms $\a:
A\>\M(A\o_{max}S)$. In this case, the exactness condition is satisfied and our theory
applies in this case too.

\begin{lema}Let $S$ be an exact $C^*$-algebra; for every $C^*$-algebra $A$ and for every
$k\geq 1$ there is a natural injective *-homomorphism $$i_A^k:\M(\U^k A\o S)\> \U^k
\M(A\o S).$$
\end{lema}
{\bf Proof.} It suffices to consider the case $k=1$, the general one follows by a similar
argument, or by an induction argument. Note first that there is a natural $S$-strict
topology on $C_b(T,A\o S)$ defined by the family of seminorms $\|f\|_s = \|f(1\o s)\| +
\|(1\o s)f\|$, for all $s\in S$. The natural inclusion map $$j_A:C_b(T,A)\o S \>
C_b(T,A\o S)$$ is continuous for this $S$-strict topology. It induces a injective
morphism between the corresponding $S$-completions
$$\bar j_A:\M(C_b(T,A)\o S) \> \overline{C_b(T,A\o S)}^{S}\simeq C_b(T,\M(A\o S)), $$
the functions on $T$ with values in $\M(A\o S)$ which are  bounded in norm and continuous
in the $S$-strict topology. Moreover, $j_A$ restricts to an isomorphism $C_0(T,A)\o S
\simeq C_0(T,A\o S)$ which in turn induces an isomorphism between the corresponding
multiplier algebras, $\M(C_0(T,A)\o S)$ and $C_0(T,\M(A\o S))$.
 Thanks to the exactness of $S$ and the lemma above, we can put these maps into a
 commutative
 diagram:
$$
\dgARROWLENGTH=0.2\dgARROWLENGTH \dgHORIZPAD=0.3\dgHORIZPAD
\begin{diagram}
\node{0}\arrow{e,t}{}\node{\M(C_0(T,A)\o
S)}\arrow{e}\arrow{s,r}{\simeq}\node{\M(C_b(T,A)\o S)} \arrow{s,r}{\bar j_A}\arrow{e}
\node{\M(\U A \o S)} \arrow{e} \arrow{s,r}{i_A}
\node{0}\\
\node{0}\arrow{e}\node{C_0(T,\M(A\o S))}\arrow{e}\node{C_b(T,\M(A\o S))}\arrow{e}\node{\U
\M(A \o S)}\arrow{e}\node{0}
\end{diagram}
$$
and deduce the existence of an injective map:
$$i_A:\M(\U A\o S)\> \U \M(A\o S).$$
Let $\f:A\>B$ be a morphism of \C s. There is a commutative diagram
$$
\begin{diagram}
\node{C_b(T,A)\o S}\arrow{s,l}{j_A}\arrow{e,t}{(\circ \f)\o \id_S}\node{C_b(T, B)\o S}
\arrow{s,r}{j_B}\\
\node{C_b(T,A\o S)}\arrow{e,t}{\circ (\f\o \id_S)}\node{C_b(T, B\o S)}
\end{diagram}
$$
from which follows the commutativity of the following diagram, expressing that the
construction of the map $i_A$ is natural:
$$
\begin{diagram}
\node{\M(\U A\o S)}\arrow{s,l}{i_A}\arrow{e,t}{(\U \f)\o \id_S}\node{\M(\U B\o S)}
\arrow{s,r}{i_B}\\
\node{\U \M(A\o S)}\arrow{e,t}{\U(\f\o \id_S)}\node{\U \M(B\o S).\ \Box}
\end{diagram}
$$

Having in mind  this  we can now say what an equivariant asymptotic morphism is:
\begin{defi} \em
An asymptotic morphism $\varphi :A\> \U^k B $ is $S$-equivariant if the following diagram
commutes:
$$
\begin{diagram}
\node{A}\arrow{e,t}{\varphi}\arrow{s,l}{\d_A}  \node{\U^k B}\arrow{se,t}{\U^k \d_B} \\
\node{\M(A \o S)}\arrow{e,t}{\varphi\o \id_S} \node{\M(\U^k B \o S)}\arrow{e,t}{i_B^k}
\node{\U^k \M(B \o S)}
\end{diagram}
$$
\end{defi}

{\bf Notation.\ }Assume that $(a_t)_{t\in T}$ and $(b_t)_{t\in T}$ are families of
elements in a \C \ $D$. We note $a_t \sim b_t$ if ${\displaystyle \lim_{t\longrightarrow
+\infty}\|a_t-b_t\|=0}$.

For an asymptotic morphism $\varphi :A \> \U B$ represented by a family of maps
$(\f_t)_{t\in T}$,  the equivariance condition reads:
$$(\varphi_t\o \id)(\d_A(a))(1\o s)\sim \d_B(\varphi_t(b))(1\o s)$$
for all $s\in S$. Note that the exactness of the algebra $S$ is necessary in defining the
asymptotic morphism $\f\o \id_S$, represented by the family of maps  $\f_t\o \id_S:
\M(A\o S)\>\M(B\o S)$.\\
\\
{\bf Examples} 1. In the case of the Hopf $C^*$-algebra associated to a group $G$ acting
on $C^*$-algebras $A$ and $B$, an asymptotic morphism $(\f_t)_t: A\>B$ is $G$-equivariant
if
$$\f_t(g(a))\sim g(\varphi_t(a))$$
for all $a\in A$ and uniformly on compacts in $g\in G$, the definition used by Guentner,
Higson and Trout in \cite{GHT}.

When the group $G$ is compact, one can average over it and find a representative for the
asymptotic morphism $\f$ such that for every $t\in T$ the map $\f_t:A\>B$ is equivariant,
in the sense that $\f_t(g(a))= g(\varphi_t(a))$. The dual
affirmation also holds as follows:\\
 2. In the case of the Hopf $C^*$-algebra associated to a coaction by
a discrete group $G$, i.e., $S= C^*_r(G)$,  an asymptotic morphism $\f :A \> \U B$ is
$S$-equivariant if $\varphi(A_g)\subset \U B_g$. Let  $(\f_t)_t: A\>B$ be any
representative of $\f$; define for each $g\in G$ a family $\psi_t^g:A_g\>B_g$ by
$\psi^g_t= P_{B_g}\circ \f_t$, where $P_{B_g}:B\>B_g$ is the projection on $B_g$. The
family $\psi^g_t$ extends to a family $(\psi_t)_t:A\>B$ which is equivalent to $(\f_t)_t:
A\>B$ and such that $\psi_t(A_g)\subset B_g$.
\\

It follows from the injectivity of $i_B$ that, if the asymptotic morphism $\varphi :A\>
\U^k B $ is $S$-equivariant, $\U^k \b$ restricts to a map
$$\d_{\varphi}: \im(\varphi) \> \M( \im (\varphi)\o S).$$
Consider now $(e_i)_i$ an approximate unit for $A$, then $\f (e_i)$ is an approximate
unit for $\im(\f)$. $\d_A(e_i)$ is an approximate unit for $\M(A\o S)$, so
$\d_{\varphi}(\f(e_i))=(\f\o \id_S)(\d_A(e_i))$ is a approximate unit for $(\f\o
\id_S)(\M(A\o S))=\M(\im \f \o S)$, hence $\d_{\f}$ is nondegenerate. Thus $\U^k \d_B$
restricts to an action on the image of an $S$-equivariant asymptotic morphism.

In the case of an action of a group $G$ this parallels the fact that the image of a
$G$-continuous element is also $G$-continuous. Hence, even if the action on $\U B$ is not
continuous in general, one can restrict to the subalgebra of $G$-continuous elements of
$\U B$. We now describe how this works in the case of a coaction of a Hopf $C^*$-algebra
$S$.

\begin{lema} Let $A$ by a $S$-algebra. The *-homomorphism $\alpha_A :A\> \U A$ is
$S$-equivariant.
\end{lema}
{\bf Proof.} In the diagram expressing the condition that $\alpha_A$ is equivariant,
$i_A$ is an isomorphism from $\M(\alpha_A(A)\o S)$ to the classes of constant functions
in $\U \M(A\o S)$ and both sides associate to an element $a\in A$ the class of the
constant function $t\mapsto \d_A(a) \in \U \M(A\o S)$. $\Box$
\begin{defi} \em
Let $A$ be an $S$-algebra, define $\U^k_S A$ the subalgebra of $\U^k A$ generated by all
elements $x$ for which there is a subalgebra  $D_x\subset \U^k A$ with the property that
$\U^k \d_B$ restricts to a nondegenerate *-homomorphism $\d_{D_x}: D_x \> \M( D_x\o S)$.
\end{defi}

It follows from the preceding lemma that the algebra $\U^k A$ is nonzero. The properties
of this construction are summarized as follows:
\begin{prop}(i)\ For every $S$-algebra $A$, $\U^k_S A$ is an  $S$-algebra with
comultiplication
given by the restriction of $\U^k \d_A$ to ${\U^k_S A}$.\\
(ii)\ $\U^k_S $ is a functor from the category of $S$-algebras and $S$-equivariant
morphisms
to itself.\\
(iii)\ There is a natural equivariant injection of the $k$-fold composition of the
functor $\U_S $ with itself into $\U^k_S $.

\end{prop}
{\bf Proof.}(i) The map $\U^k \d_A$ which restricted to each $D_x$ gives an $S$-algebra
structure extends to the $C^*$-algebra generated by all this subalgebras of $\U^k A$. To
see that this *-homomorphism is nondegenerate observe that $[\U^k \d_A(D_x)]( D_x \o S)$
is dense in $D_x \o S$ hence $[\U^k \d_A (\U^k_S A)](\U^k_S A\o S)$ is dense in $D_x \o
S$ and thus also in $\U_S^k A\o S$.

For (ii) we restrict to the case $k=1$, the general one works  in a similar way. Let
$D\subset \U A$  be an algebra for which  $\U \d_A $ restricts to an action of $S$, thus
in the following diagram:
$$
\begin{diagram}
\node{D}\arrow{e}\arrow{s}\node{\U A}\arrow{e,t}{\U \varphi}\arrow{s,l}{\U \a}
\node{\U B}\arrow{s,r}{\U \d_{B}} \\
\node{\M(D\o S)}\arrow{e}\arrow{se}\node{\U \M(A \o S)}\arrow{e,t}{\U(\varphi \o \id_S)}
\node{\U \M( B \o S)}\\
\node{} \node{\M(\U A \o S)}\arrow{n,l}{i_A} \arrow{e,t}{\U \f\o \id_S } \node{\M(\U B\o
S)} \arrow{n,r}{i_B}
\end{diagram}
$$
the left side commutes. The upper right square is commutative because it is the image by
$\U$ (which is a functor on the category of \C s) of the square expressing the
$S$-equivariance of $\f$. Finally, the right lower square commutes by the naturality of
the applications $i_A$ and $i_B$. It follows that $\U\f(D)\subset \U_S B$, which implies
that $\U\f(\U_S A)\subset
\U_S B$, and also that $\U \f$ is $S$-equivariant.\\
(iii) It is enough to consider the case $k=2$; let $D$ be a subalgebra of $\U(\U_S A)$
which are among those generating $\U_S(\U_S A)$, there is a commutative diagram:
$$
\begin{diagram}
\node{D}\arrow{e}\arrow{s}\node{\U(\U_S A)}\arrow{e,t}{}\arrow{s,l}{}
\node{\U^2 A}\arrow{s}\\
\node{\M(D\o S)}\arrow{e}\arrow{e}\node{\U \M(\U_S A \o S)}\arrow{e,t}{} \node{\U^2 \M( A
\o S)}
\end{diagram}
$$
with horizontal maps given by inclusion and vertical maps by the restriction of $\U^2
\d_A$. It follows that $D\subset \U^2_S A$ and, by the definition of $\U_S$, that
$\U_S(\U_S A)\subset \U^2_S A$. The equivariance and the naturality of this inclusion
follow from the corresponding commutative diagrams.
$\Box$\\

Given an $S$-algebra $A$, denote by $\Sigma A$, the suspension of $A$, i.e. the tensor
product algebra $C_0(\mathbb{R})\o A $, with the trivial coaction of $S$ on
$C_0(\mathbb{R})$.

One key ingredient of $E$-theory is the Connes-Higson construction of an asymptotic
morphism  which associates to  a short exact sequence $0\>J\>B\>A\>0$ of separable
$C^*$-algebras an asymptotic morphism $\partial : \Sigma A\> \U J$ from the suspension of
$A$ to $J$.

We now show that the same construction exists also in the $S$-equivariant setting.

The Connes-Higson asymptotic morphism is based on the following notion:

\begin{defi} \em Let $J$ be an ideal in a $C^*$-algebra $B$; a quasicentral
approximate unit for $J\vartriangleleft B$ is a norm continuous family  $(u_t)_{t\in T}$
of elements of
$J$ such that\\
(i) $0\leq u_t \leq 1$,\\
(ii) $u_t j \sim j$,\\
(iii) $\ u_t b \sim b u_t$,\\
for all $j\in J$ and $b\in B$.
\end{defi}

The Connes-Higson asymptotic morphism is defined as the the asymptotic family which, for
$f\in C_0(0,1)$ and $a\in A$, is given by
$$\partial_t(f\o a)= f(u_t)q(a)$$
for some choices of an approximate unit $(u_t)_{t\in T}$ for $J\vartriangleleft B$ and of
a section $q:A\>B$. It turns out that the homotopy class of this asymptotic morphism do
not depend on these choices so the formula defines an element of $\{\Sigma A,J\}$.

We assume now that $B$ is an $S$-algebra and that the coaction of $S$ on $B$ restricts to
a coaction on the ideal $J$.

\begin{defi} \em An $S$-quasicentral approximate unit for $J\vartriangleleft B$ is a
quasicentral approximate unit $(u_t)_{t\in T}$ such that $$\d (u_t)(1\o s)\sim (u_t\o
1)(1\o s)$$ for all $j\in J$, $b\in B$ and $s\in S$.
\end{defi}

It was proven in \cite[Lemma 4.1]{BS} that such approximate units exist provided that $B$
and $S$ are separable. In the case of an action of a group $G$, this latest condition
states that $g(u_t)\sim u_t$, uniformly on the compacts of $G$.

\begin{teo}Let $0\>J\>B\>A\>0$ be an $S$-extension; the Connes-Higson asymptotic
morphism $\partial : \Sigma A\> \U J$ associated to it is $S$-equivariant.
\end{teo}
{\bf Proof.} We have to verify the equivariance condition which in this case is :
$$i_S\circ(\partial \o \id_S) \circ \d_{\Sigma A} = \U \d_J \circ \partial$$
in $\U \M(J\o S)$ or in an equivalently that $\d_J(\partial _t(f\o a))(1\o s) \sim
(\partial_t \o \id) (\d_{\Sigma A}(f\o a))(1\o s)$ for every $f\in C_0(0,1)$, $a\in A$
and $s\in S$.

The action on $\Sigma$ is trivial thus
$$(\partial_t\o \id_S)(\d_{\Sigma A}(f\o a)) = (\partial_t\o \id_S )(f\o\d_A(a)).$$
The family of maps $\partial_t\o \id_S:\Sigma A\o S\> J\o S$ is asymptotically equivalent
to the family associated to the short exact sequence:
$$0\>J\o S\> B\o S\> A\o S\>0$$
\cite[Proposition 5.9]{GHT}, which in turn is given by $f(u_t\o v_t)\tilde q(\d_A(a))$,
where $(v_t)_{t\in T}$ is an approximate unit for the $C^*$-algebra $S$, and $\tilde q:
A\o S\> B\o S$ is a section.

Note that
$$f(u_t\o v_t){\tilde q}(\d_A(a))\sim f(u_t\o v_t)\d_B(q(a))$$
as $\tilde q(\d_A(a))-\d_B(q(a))\in \ker (p\o \id_S)=J\o S$ and $f(u_t\o v_t)h\sim 0$
 for all $h\in
J\o S$ (see \cite[Lemma 5.6]{GHT}). $f(u_t\o v_t)\sim f(u_t\o 1)$ in the $S$-strict
topology of $\M(J\o S),$ as the condition $\|[f(u_t\o v_t)-f(u_t\o 1)](1\o s)\|\> 0$ is
easily verified for $f(x)=x$ and by the Weierstrass approximation theorem holds for all
$f\in C_0(0,1)$. Consequently
$$(\partial_t\o \id_S)(\d_{\Sigma A}(f\o a))(1\o s)\sim f(u_t\o 1)\d_B(q(a))(1\o s),$$
where we regard $\d_B(q(a))\in\M(B\o S)$ as an element of $M(B\o S)$ and, by restriction,
as an element of $M(J\o S)$.

On the other hand, $\d_J(\partial_t(f\o a))=f(\d_J(u_t))\d_B(q(a))$, hence we have to
prove that $f(\d_J(u_t))\sim f(u_t\o 1)$ in the strict topology of $\M(J\o S)$ which
follows again from Weierstrass's theorem and the $S$-equivariance of $(u_t)_{t\in T}$.
$\Box$\\

Like in the non-equivariant case, the homotopy class of the Connes-Higson asymptotic
morphism does not depend on the choices of the approximate unit $u_t$ and of the section
$q$.

It is worth mentioning that, up to suspension, every asymptotic morphism is obtained in
this way from an extension, hence the theorem shows that the definition of an
$S$-asymptotic morphism is appropriate.

\section{E-theory}

In the nonequivariant case, $E$-theory groups are defined as
$$E(A,B)=\{\Sigma A\o \K,\Sigma B\o \K\},$$
the homotopy classes of asymptotic morphisms between the algebras $\Sigma A\o \K$ and
$\Sigma B\o \K$. In the case of a coaction by a Hopf algebra $S$, a similar definition of
$E^S$, using $S$-asymptotic morphisms instead of nonequivariant ones, would have many
good properties, but not a crucial one: to insure the existence of a natural
transformation, $KK^S\>E^S$, from equivariant $KK$-theory to it. In order to have such a
property, one has to define $E^S$ in such a way that it has the same stability properties
as $KK^S$.

In the case of actions by a group $G$, one can take advantage of the the $G$-module
$L^2(G)$ endowed with the left regular representation.  It has the following property (
\cite[Lemma 2.3]{MP} ): let $A$ be a $G$-algebra, and let $E_1$ and $E_2$ be Hilbert
$G$-$A$-modules which are (non-equivariantly) isomorphic as Hilbert $A$-modules, then
$L^2(G,E_1)$ and $L^2(G,E_2)$ are isomorphic Hilbert $G$-$A$-modules.

Here a module $L^2(G,E)$ denotes the tensor product of the modules $L^2(G)$ and $E$.

Denote by $\K_G$ the $G$-algebra of compact operators of the module $L^2(G)^{\infty}$,
with the action of $G$ induced from the action on $L^2(G)$, denoted by $\lambda$.

In the case of an action of a group $G$, the equivariant $E$-theory, $E^G$, is defined in
\cite{GHT} as
 $$E^G(A,B)=\{\Sigma A\o \K_G,\Sigma B\o \K_G\}.$$
As explained in the appendix, this insures equivariant stability.

In the more general case of a coaction by a Hopf $C^*$-algebra $S$, we do not have a
substitute for the module $L^2(G)$.

To solve this problem, we construct a category of $S$-algebras in which we allow a bit
more freedom for the coaction on the algebras $A\o\K$ besides the coaction $\a\o
\id_{\K}$, in the sense of the following definition. We shall show in the next section
that this is exactly what we need.

\begin{defi} \em Let $(S,\d_S)$ be a Hopf-$C^*$-algebra, and $(A,\a)$ be an $(S,\d_S)$-algebra.
 A coaction $\d_{A\o \K}:
A\o\K\>\M(A\o\K\o S)$ of $S$ on $A\o\K$ is compatible with the coaction $\a$ of $S$ on
$A$ if there is a minimal projection $e\in \K$ such that the map $a\>a\o e$ is
equivariant.
\end{defi}

A key notion is the following:

\begin{defi} \em Let $(S,\d_S)$ be a Hopf-$C^*$-algebra, and $(A,\a)$ be an $(S,\d_S)$-algebra.
A $\a$-cocycle is an unitary $V\in M(A\o S)$ such that $$(\id_A\o\d_S)V=(V\o \id_S)[(\a\o
\id_S)(V)].$$
\end{defi}

{\bf Example.}  Let $G$ be a locally compact group and $C_0(G)$ the associated Hopf
algebra. A $\a$-cocyle is a strongly continuous map from $G$ into the unitary group of
$A$ such that $V_{g_1g_2}=V_{g_1}\alpha_{g_1}(V_{g_2})$. (cf. \cite[Def. 2.2.3
]{Co|classiffactIII})

\begin{lema} If $V$ is a $\a$-cocyle, then $\a^V(\cdot):= V\a(\cdot)V^*$
is a coaction of $S$ on $A$.
\end{lema}

{\bf Proof.} The $S$-strict topology is stronger than the strict topology induced from
$M(A\o S)$, so comparing completions yields an identification $M(\M(A\o S))=M(A\o S)$.
Hence $\a^V(\cdot): A\>\M(A\o \K)$ is a nondegenerate map. The equality $(\id_A \o
\d_S)\a^V(\cdot)= (\a^V\o \id_S)\a^V(\cdot)$ is easy to check. $\square$
%
%\begin{eqnarray}
%(\id_A \o \d_S)\a^V(\cdot) &=& (\id_A \o \d_S)(V\a (\cdot)V^*)  \nonumber \\
%&=& (\id_A \o \d_S)(V)(\id_A \o \d_S)\a(\cdot)(\id_A \o \d_S)(V)^* \nonumber \\
%&=& (V\o \id_S)(\a\o \id_S)(V)(\id \o \d_S)\a(\cdot) [(V\o \id_S)(\a\o \id_S)(V)]^* \nonumber \\
%&=& (V\o \id_S)(\a\o \id_S)(V)(\a\o \id_S)\a(\cdot) (\a\o \id_S)(V)^*(V\o \id_S)^* \nonumber \\
%&=& (V\o \id_S)(\a\o \id_S)(V\a(\cdot)V^*)(V\o \id_S)^* \nonumber \\
%&=& (V\a(\cdot)V^*\o \id_S)(V\a(\cdot)V^*) \nonumber \\
%&=& (\a^V\o \id_S)\a^V(\cdot).\ \square \nonumber
%\end{eqnarray}

Note that if $\f:(A,\a)\> (B,\b)$ is a surjective $S$-morphism and  $V$ is a
$\a$-cocycle, then $(\f\o \id_S)(V)$ is a $\b$-cocycle.

The following lemma is a generalization of a result due to Connes in the case of a group
acting on a $C^*$-algebra.

\begin{lema} Let $(A, \d_A)$ be an $S$-algebra; there is a coaction
$\d_{M_2(A)}:M_2(A)\>\M(M_2(A)\o S)$ of $S$ on $M_2(A)$ with the property that for all
$a,b\in A$,
\begin{displaymath}
\d_{M_2(A)} \left(\begin{array}{cc}
a & 0\\
0 & b
\end{array}\right)=
\left(\begin{array}{cc}
\a(a) & 0\\
0 & \a'(b)
\end{array}\right),
\end{displaymath}
for some coaction $\a'$ of $S$ on $A$,  if and only if there is a $\a$-cocycle
 $V \in M(A\o S)$ such that $\d_{M_2(A)}$ is given by
 $\d_{M_2(A)}\ma a b c d = \ma {\a(a)} {\a(b)V^*} {V\a(c)} {V\a(d)V^*}$.
 \end{lema}
{\bf Proof.} The converse is obvious. Let $\ma A B V D$ be $\d_{M_2(A)}\ma 0 0 1 0$. It
follows easily that $A=B=D=0$, and that $V$ is an unitary in $M(A\o S)$.
\begin{eqnarray}
(\id_{M_2(A)}\o \d_S)\ma 0 0 V 0 &=& (\id_{M_2(A)}\o \d_S)\d_{M_2(A)}\ma 0 0 1 0 \nonumber \\
 &=& (\d_{M_2(A)}\o \id_S)\d_{M_2(A)}\ma 0 0 1 0\nonumber \\
  &=& (\d_{M_2(A)}\o \id_S)\ma 0 0 {1_{A\o S}} 0(\d_{M_2(A)}\o \id_S)\ma V 0 0 0\nonumber \\
   &=& \ma 0 0 V 0 \o \id_S \ma {(\a\o \id_S) V} 0 0 0\nonumber \\
   &=& \ma 0 0 {(V\o \id_S)(\a \o \id_S)V } 0,\nonumber
\end{eqnarray}
 which proves that $V$ is an $\a$-cocycle. $\square$
\\
The coactions $\d_A$ and $\d_A'$ like in the lemma above are called {\it exterior
equivalent}.\\

An advantage of working with cocycles  is that, given a coaction $\a$ of $S$ on $A$, we
can control a compatible coaction on $A\o \K$. This goes as follows: let $ \d_{A\o\K}:
A\o\K\>\M(A\o\K\o S)$ be such a compatible action and let $\{e_{ij}|\ i,j\geq1\}$ be a
set of matrix units for $\K$ arranged such that $e_{11}=e$; define
$$V= \sum (1_A\o e_{i1}\o 1_S)\d_{A\o\K}(1\o e_{1i}).$$

\begin{lema} $V$ is a $\d_{A\o\K}$-cocycle and
$V\d_{A\o\K}(\cdot)V^*=\a(\cdot)\o \id_{\K}$.
\end{lema}
{\bf Proof.} It is easy to see that $V$ is an unitary; let us check the cocycle
condition:
\begin{eqnarray}
(V\o \id_S)[(\d_{A\o \K}\o \id_S)(V)] &=& \sum_{i,j}[1_A\o e_{i1}\o 1_{S\o S}][\d_{A\o
\K}(1_A\o e_{1i})\o 1_S]
 [\d_{A\o \K}(1_A\o e_{j1})\o 1_S] \nonumber \\
 & & [(\d_{A\o \K}\o \id_S)\d_{A\o \K}(1\o e_{j1})]   \nonumber \\
 &=&  \nonumber \sum_{i,j}[1_A\o e_{i1}\o 1_{S\o S}][\ab((1_A\o e_{1i})(1_A\o e_{j1}))\o 1_S]\\
 & & [(\id_{A\o\K}\o \d_S)\ab(1\o e_{j1})] \nonumber\\
  &=& \sum_i[1_A\o e_{i1}\o 1_{S\o S}][\ab(1_A\o e_{11})\o 1_S]\nonumber \\
 & & [(\id_{A\o\K}\o \d_S)\ab(1\o e_{i1})]\nonumber \\
   &=& \sum_i[1_A\o e_{i1}\o 1_{S\o S}][1_{A\o S}\o e_{11}\o 1_{S\o S}]
 [(\id_{A\o\K}\o \d_S)\ab(1\o e_{i1})]\nonumber \\
 &=& \sum_i[1_A\o e_{i1}\o 1_{S\o S}][1_A\o e_{11}\o 1_S]
 [(\id_{A\o\K}\o \d_S)\ab(1\o e_{i1})]\nonumber \\
 &=& (\id_{A\o\K}\o \d_S)(\sum_i(1_A\o e_{i1}\o 1_{S})\ab(1_A\o e_{i1}\o 1_{S\o S}))\nonumber \\
 &=& (\id_{A\o\K}\o \d_S)V. \nonumber
\end{eqnarray}
Also, for every $p,q\geq 1$:
\begin{eqnarray}
V\ab(a\o e_{pq})V^* &=&  \sum_{i,j}(1_A \o e_{i1}\o 1_S)\ab(1_A \o e_{1i})
 \ab(1_A \o e_{pq})\ab(1_A \o e_{j1}) \nonumber \\
&& (1_A \o e_{1j}\o 1_S) \nonumber \\
&=& (1_A \o e_{p1}\o 1_S)\ab(1_A \o e_{11})(1_A \o e_{1q}\o 1_S) \nonumber \\
&=& (1_A \o e_{p1}\o 1_S)(\a(1_A) \o e_{11})(1_A \o e_{1q}\o 1_S) \nonumber \\
&=& (\a\o \id_{\K})(a\o e_{pq}), \nonumber
\end{eqnarray}

which proves that $V\a(\cdot)V^*=\a(\cdot)\o \id_{\K}$. $\square$

\begin{defi} \em Let $(A,\a)$ and $(B,\b)$ be $S$-algebras, we call $\bar S$-morphism an
$S$-equivariant morphism  $\f:(A\o\K, \d_A\o\id_\K)\>(B\o\K, \d_{B\o \K})$, where
$\d_{B\o \K}$ is a coaction of $S$ on $B\o\K$ which is compatible with $\d_B$.
\end{defi}

For an interval $I=[a,b]$ and an $S$-algebra $A$, denote by
$$IA=\{f:I\>A \ |\ f \ \mathrm{ is\ continuous}\},$$
the algebra of continuous functions from $I$ to $A$, regraded as an $S$-algebra with the
trivial coaction of $S$ on $I$.

\begin{defi} \em Let $(A, \d_A)$ and $(B, \d_B)$ be $S$-algebras, two $\bar S$-morphisms
$\f_0, \f_1: (A\o \K, \d_A\o\id_\K)\> (B\o\K,\d_{B\o \K})$ are $S$-homotopic, if there is
an $S$-morphism $\Phi:(A\o\K,\d_A\o\id_\K) \>(I(B\o\K),\id_I\o \d_{B\o\K})$ such that
$\f_0$ and $\f_1$ are obtained by evaluation of $\Phi$ at the endpoints of the interval
$I$.
\end{defi}

We shall write $\f_0\sim_{\bar S}\f_1$ in this case; $S$-homotopy is an equivalence
relation.

Recall that a morphism $\f:A\>B$ is {\it quasiunital} if there is an approximate unit
$(e_i)_i$ for $A$ such that $(\f(e_i))_i$ converges strictly to a projection $p\in M(B)$.
This condition is equivalent to the existence of an extension of $\f$ to a strictly
continuous map $\bar\f:M(A)\>M(B)$ between the multiplier algebras.

We say that an $S$-algebra $(A,\a)$ is {\it trivially stable}, if there is an
$S$-equivariant isomorphism $(A, \d_A)\simeq (A\o \K, \d_A\o \id_\K)$.

\begin{lema} Let $(A,\d_A)$, $(B, \d_B)$ be $S$-algebras, and assume that $B$
is  trivially stable,
then every $\bar S$-morphism between $A$ and $B$ is $\bar S$-homotopic to a quasiunital
one.
\end{lema}
{\bf Proof.} In the non equivariant case this is  \cite[Theorem 1.3.16]{JT}. In our case,
take $\f:A\o\K\>B\o\K$ to be a $*$-homomorphism, $\a\o \id_\K $-$\bb$-equivariant for
some coaction $\d_{B\o \K}$ on $B\o\K$. Because $B$ is trivially stable and $\d_{B\o \K}$
is compatible with $\d_B$, it follows that $$(B\o\K,\bb)\simeq ((B\o\K)\o \K,\d_{B\o
\K}\o \id_\K).$$  One can check now that the decompositions in the proof of the theorem
just quoted carry over to our case, as the action on the last copy of $\K$ is trivial.
$\square$

{\it From now on, we assume in this section that all algebras are trivially stable.} This
is not a restriction as, for our purposes, one can replace any algebra by its tensor
product by the compacts, and have a trivially stable one. Thus when we write $A\o \K$
with some compatible coaction $\d_{A\o \K}$, we imply that the algebra $A$ is also
trivially stable, but we do not specify an extra copy of $\K$.

\begin{prop} There is  a category, denoted $\ss$, whose objects are $S$-algebras
and whose morphisms are $S$-homotopy classes of $\bar S$-morphisms.
\end{prop}
{\bf Proof.} First suppose that  $\f:A\o\K\>B\o\K$ is a quasi-unital $*$-homomorphism
which is $\ab$-$\bb$-equivariant, for a coaction $\ab$ on $A\o\K$ compatible with the
given coaction $\a$ on $A$, and a coaction $\bb$ on $B\o\K$, compatible with the given
coaction $\b$ on $B$. Then $\f$ is among our $\bar S$-morphisms as follows.

Denote by $p:=\f(1)$; $p$ is a projection in the the center of the multipliers $M(B\o
\K)$ and $B\o\K$ decomposes as a direct sum of algebras $p(B\o\K)\oplus(1-p)(B\o\K)$.
Moreover, $p$ and $1-p$ are $\bb$-equivariant and there is an identification
$$(p\o 1_S)\M(B\o\K\o S)=\M(p(B\o\K)\o S)$$
which is the $S$-completion of the identification $(p\o 1_S)(B\o \K\o S)=p(B\o \K)\o S$.
It follows that the action $\bb$ decomposes as a sum
$$p\bb\oplus (1-p)\bb: p(B\o \K)\oplus
(1-p)(B\o \K)\>\M(p(B\o \K)\o S)\oplus\M((1-p)(B\o \K)\o S).$$
 Denote by $U\in M(B\o\K \o S)$
the $\ab$-cocycle for which $\ab(\cdot)=U(\a\o \id_\K)(\cdot)U^*$; using the
decomposition above we denote by $\f^\sharp(U)$ the unitary $(\f\o \id_S)(U)\oplus 1$. It
follows that $\f^\sharp(U)$ is a $\bb$-cocycle and that  $\f$ is $\a\o
\id_\K(\cdot)$-$\f^\sharp(U)\bb(\cdot)\f^\sharp(U)$-equivariant, which proves the claim.

For the composition, let $\f:(A,\a\o \id_{\K})\>(B,\d_{B\o\K})$ and $\psi:(B,\b\o
\id_\K)\>(C,\d_{C\o\K})$ be two $\bar S$-morphisms. We can write $\bb(\cdot)=U_B(\b\o
\id_K)(\cdot)U_B^*$ for some $\bb$-cocycle $U_B$. If follows that the morphism
$\psi\circ\f$ is $(\a\o \id_\K)(\cdot)$- $\psi^\sharp(U_B)(\c\o
\id_\K)(\cdot)(\psi^\sharp(U_B))^*$-equivariant, which concludes the proof. $\square$

Denote by $\{A,B\}^{\bar S}_0$ the homotopy classes of $\bar S$-morphisms between $A$ and
$B$, i.e., the morphisms in $\ss$ between $A$ and $B$.

\begin{defi} \em Let $A$, $B$ be two $S$-algebras. An $\bar S$-asymptotic morphism is an
$S$-morphism $\f:A\>\U^n_S B$, for some $n \geq 1$.
\end{defi}

Let $(A, \d_A)$ and $(B, \d_B)$ be $S$-algebras and take two $\bar S$-asymptotic
morphisms between them, $\f_0$ and  $\f_1$; we say that they are $n$-$S$-homotopic, if
there is an $S$-morphism $\Phi:(A\o\K,\d_A\o\id_\K) \>\U^n_S((I(B\o\K)),\id_I\o
\d_{\U^n_S(B\o \K)})$ such that $\f_0$ and $\f_1$ are obtained by evaluation of $\Phi$ at
the endpoints of the interval $I$, for some compatible coaction $\d_{\U^n_S(B\o \K)}$ on
$B\o \K$.

It follows like in the group case, \cite[Proposition 2.3]{GHT}, that
\begin{lema}
$n$-homotopy is an equivalence relation on  $\bar S$-asymptotic morphisms from $A$ to
$\U^n_SB$.
\end{lema}

Denote by $\{A,B\}_n^{\bar S}$ the set of $n$-homotopy classes of $\bar S$-asymptotic
morphisms from $A$ to $\U^n_SB$.

Given  an $\bar S$-asymptotic morphism $\varphi:A\longrightarrow\U^n_S B$,  the
composition
$$A\stackrel{\varphi}{\longrightarrow}\U^n_S B
\stackrel{\U^n(\alpha_B)}{\longrightarrow}\U^{n+1}_S B$$ is an $\bar S$-asymptotic
morphism; this map agrees with $n$-$S$-homotopy of asymptotic morphisms and thus induces
a map $\{A,B\}_n^{\bar S}\> \{A,B\}_{n+1}^{\bar S}$; an equivalent way of defining this
map is
 $A\stackrel{\varphi}{\longrightarrow}\U^n_S B
\stackrel{\alpha_{\U^nB}}{\longrightarrow}\U^{n+1}_SB.$
\begin{defi} \em
Denote by $\{A,B\}^{\bar S}$ the inductive limit ${\displaystyle \lim_{\longrightarrow}}
\{A,B\}_n^{\bar S}$.
\end{defi}
\begin{teo}
Let $A,B$, and $C$ be $S$-algebras. Given $\bar S$-asymptotic morphisms
$\varphi:A\>\U^j_S B$ and $\psi:B\>\U^k_S C$, the formula
$A\stackrel{\varphi}{\longrightarrow}\U^j_S B
\stackrel{\U^j(\psi)}{\longrightarrow}\U^{j+k}_S C$,  defines an  associative composition
law
$$\{A,B\}^{\bar S}\times\{B,C\}^{\bar S}\longrightarrow \{A,C\}^{\bar S}.$$
Moreover, if $A$ is separable, the inclusion $\{A,B\}_1^{\bar S}\>\{A,B\}^{\bar S}$ is a
bijection.
\end{teo}
{\bf Proof.} The proof of the first part is exactly the same as in the group case
\cite[Proposition 2.12]{GHT}: the homotopies used in the proof can be used in the
$S$-equivariant case.

The second part is an alternative formulation of the composition of asymptotic morphisms
given in \cite{CH}, and adapts as follows.

The key ingredient is the following result: let $D\subset \U^2 C$ be a separable algebra,
there exists a function $r:T\>T$ such that the restriction of a two variables function to
the graph of $r$ defines a *-homomorphism $R:D\> \U C$. Moreover, the inclusion $D\subset
\U^2 C$ is 2-homotopic to the composition $\alpha_C\circ R$ and hence factorizes through
$\U C$. This is the first part of \cite[Lemma 2.17]{GHT}. The second part states that
this choice of the function $r$ can be done with respect to a map between separable
subalgebras $D_1\>D_2$. The $S$-equivariant case follows from this by taking
$D_2=\d_{D_1}(D_1)\subset\M(D_1\o S)$.

The $2$-homotopy used in this factorization is given by
$$H(t_1,t_2,s)= \left\{\begin{array}{ccc}
F(t_1,t_2)&\mbox{if}& t_1>sr(t_2) \\
F(sr(t_2),t_2)&\mbox{if}& t_1\leq sr(t_2)
\end{array}\right. $$
and it is an $S$-homotopy. $\square$

We can now define the $E$-theory groups.

\begin{defi} \em Let $A$ and $B$ be $S$-algebras; we define the equivariant $E$-theory
of $A$ and $B$ as
$$E^S(A,B)=\{\Sigma A,\Sigma B \}^{\bar S}.$$
\end{defi}

An asymptotic morphism $\f_t:\Sigma A\o \K \>\Sigma B\o \K$ for which there is a norm
continuous family of $\d_B\o \id_{\Sigma \o \K}$-cocycles $U_t\in M(\Sigma B\o \K \o S)$
such that
$$(\f_t\o \id_S)(\a(a))\sim U_t\b(\f_t(a))U_t^*$$
in the $S$-strict topology will define an element of $E^S(A,B)$.
\\
{\bf Example.} Let us take a closer look at the case the group case; assume that the
$C^*$-algebras $A$ and $B$ are endowed with coaction of a group $G$, denoted respectively
with $\alpha,$ and $\beta$.

Let  $\f_t:A\o\K\>B\o\K$ be an asymptotic morphism, and assume that
$$\f_t(\alpha_g(a))\sim \beta^t_g(\f_t(a)),$$ where  $\beta^t$ is a continuous
family of actions on $B\o \K$, all exterior equivalent with the given action $\beta$.
Such an element defines an element of the $E$-theory group $E^{C_0(G)}(A,B)$.

An important example of such an asymptotic morphism appears in the proof of the
Baum-Connes for amenable groups where the family $\beta^t$ is due to a deformation of the
action on a Hilbert space $H$, reflected on the algebra of compact operators $\K=\K(H)$ (
see \cite[Remark 2.7.6, Proposition 4.6.25]{GH|LNM}).
\\
\\
As in the nonequivariant case, $E^S(A,B)$ is an abelian group. The sum is defined using
an isomorphism $M_2(\K)\simeq \K$ and the inverse is defined using the fact that an
element of the type $t\mapsto\ma {f(t)} 0 0 {f(-t)}$ from $C_0(\mathbb{R})$ to a
$C^*$-algebra is null homotopic.

From Theorem 3.13 follows that:

\begin{prop} For any separable $S$-algebras $A$, $B$ and $C$, there is a bilinear
composition law
$$E^S(A,B)\times E^S(B,C)\>E^S(A,C)$$
extending the composition of the category $\ss$. $\square$
\end{prop}

Any morphism of Hopf-\C s $\f:S\> M(S')$ induces a restriction natural transformation
$E^S(A,B)\> E^{S'}(A,B)$.

It follows form the discussion the the appendix that for every $S$-algebra $A$, the
functors $E^S(A, \cdot)$ and $E^S(\cdot, A)$ are stable, in the equivariant setting.

$E$-theory appeared as an answer to the question whether there are always six-terms exact
sequences in $KK$-theory. The key property in establishing this is the following notion:
a functor $F$ from a category of algebras to abelian groups is {\it half-exact} if given
a short exact sequence
$$0\>J\>B\>A\>0$$ of algebras, the sequence $F(J)\>F(B)\>F(A)$ is exact.
We saw that the Connes-Higson asymptotic morphism is $S$-equivariant, the proof from the
nonequivariant case adapts to show that given an $S$-algebra $D$, the functors
$E^S(D,\cdot)$ and $E^S(\cdot,D)$ are half-exact. Applying the Puppe exact sequences
techniques along with Cuntz' theorem on Bott periodicity yields:

\begin{prop}
Let $D$ be an  $S$-algebra and let
$$0\>J\>B\>A\>0$$
be an $S$-extension. There are six-terms exact sequences:
$$
\begin{diagram}
\node{E^S(D,SJ)}\arrow{e} \node{E^S(D,SB)}\arrow{e}
\node{E^S(D,SA)}\arrow{s}\\
\node{E^S(D,A)}\arrow{n} \node{E^S(D,B)}\arrow{w} \node{E^S(D,J)}\arrow{w}
\end{diagram}
$$
and
$$
\begin{diagram}
\node{E^S(SJ,D)}\arrow{s} \node{E^S(SB,D)}\arrow{w}
\node{E^S(SA,D)}\arrow{w}\\
\node{E^S(A,D)}\arrow{e} \node{E^S(B,D)}\arrow{e} \node{E^S(J,D). \square}\arrow{n}
\end{diagram}
$$
\end{prop}

\begin{prop} For any exact  $S$-algebra $D$ there is a natural transformation
$$\tau_D : E^S(A,B)\> E^S(A\o D,B\o D)$$
such that $\tau_D(1_A)= \tau_{A\o D}$.
\end{prop}

{\bf Proof.} It is shown in \cite[Proposition 4.4]{GHT}, that, thanks to the exactness of
$D$, there is a natural map $i_D: \U^k(A\o D)\>\U^k A\o D$ for all $k\geq 0$. The
argument is similar to the one used in the proof of the Lemma 2.5 and thus it is easy to
see that this map restricts to a map
$$i_D^S: \U^k_S(A\o D)\>\U^k_S A\o D.$$
Composition with $i_D^S$ gives a well defined map
$$\tau_D:\{A,B\}^{\bar S}\> \{ A\o D,B\o D \}^{\bar S},$$
because the tensor product of two compatible coactions is a compatible coaction. and
taking suspensions one obtains a map $\tau_D : E^S(A,B)\> E^S(A\o D,B\o D)$. The second
condition is obvious by the definition of $i_D$. $\square$ \\
\begin{rmq}\em Based on this, one we can define
$$\o:E^S(A,B)\times E^S(C,D)\>E^S(A\o C, B\o D)$$ by $\f\o\psi:=(\f\o1_C)\circ(1_B\o
\psi)$ for $\f\in E^S(A,B)$ and $\psi\in E^S(C,D)$.
\end{rmq}

\begin{prop}
For every group $G$ and every $G$-algebras $A, B$, the $E$-theory of $A$ and $B$ with
respect to coactions of $C_0(G)$, $E^{C_0(G)}(A,B)$, is isomorphic to the $E$-theory with
respect to actions of $G$, $E^G(A,B)$.
\end{prop}
{\bf Proof.} We already observed that $C_0(G)$-asymptotic morphisms correspond to
$G$-asymptotic morphisms  as defined in \cite{GHT}. Combining this with  the stability of
$E^{C_0(G)}$, allows us to replace in the definition of $E^{C_0(G)}(A,B)$, $(A, \alpha)$
with $(A\o \K_G, \alpha \o \lambda_G)$, and $(\U_S B\o \K, \d_{\U_SB \o \K})$ with $(\U_S
B\o \K \o \K_G, \d_{\U_SB \o \K}\o \lambda_G)$. As proven in the appendix, this latest
algebra is $G$-isomorphic to $(\U_S B\o \K_G, \U \d_{B}\o \lambda_G)$, from which the
result follows. $\square$
\\
\\
We now prove the Baaj-Skandalis isomorphism in $E$-theory. First let us remind some
definitions and properties of the cross products by actions and by coactions of groups on
$C^*$-algebras.

Let $A$ be a $C^*$-algebra endowed with a coaction $\d_A : A \> \M(A\o C^*_r(G)) $ of a
group $G$; to construct the cross product algebra with make the assumption that
$\d_A(A)(A\o C^*_r(G))$ is dense in $A\o C^*_r(G)$, and we say that $A$ is a $\hat
G$-algebra.

This condition is always fulfilled for coactions of discrete and of amenable groups.

\begin{defi} \em Let $(A,\d_A)$ be a $\hat G$-algebra; the reduced cross product by the coaction
$\d_A$ is the subalgebra of $M(A \o \K(L^2(G)))$ spanned by the products $\d_A(a)(\id_A\o
M_f)$, for all $a\in A$ and $f\in C_0(G)$, $M_f$ denoting the multiplication operator by
$f$.
\end{defi}

We denote by $A \r \hat G$ this algebra. There are inclusion maps $i_A :A \> \M(A\r \hat
G)$ given by $i_A(a):= \d_A(a)$, and $i_G : C_0(G)\> \M(A\r \hat G)$ given by $i_G(f):=
\id_A\o M_f$; the cross product is the closed linear span of the products $i_A(a)i_G(f)$
for $a\in A$ and $f\in C_0(G)$.

The algebra $A \r \hat G$ is endowed with an action of $G$, given by
$\alpha_g(i_A(a)i_G(f))=i_A(a)i_G(gf)$, where $g(f)(h):=f(gh)$.

The cross product by a coaction is a functor from the category of $\hat G$-algebras to
the category of $G$-algebras: if $A$ and $B$ are $\hat G$-algebras and if $f:A\>B$ is an
equivariant $*$-homomorphism, then there is an induced morphism $\f \times 1:A \r \hat
G\> B \r \hat G$, which is given by $\f \times 1(i_A(a)i_G(f))=i_B(\f(a))i_G(f)$.

We can now state the analogue of the Takesaki-Takai duality for noncommutative groups.

Let $A$ be a $G$-algebra, the cross product $A\r G$ is a $\hat G$-algebra, and there is a
$G$-equivariant isomorphism between $A\r G \r \hat G$ and $A\o \K(L^2(G))$.

Similarly, if $(A, \d_A)$ is a $\hat G$-algebra, there is a $\hat G$-equivariant
isomorphism between the double cross product $A\r \hat G \r G$ and $A\o \K(L^2(G))$. Here
the coaction of $G$ on $A\o \K(L^2(G))$ is given by $(\id_A\o W^*)(\id_A \o
\sigma)(\delta_A \o \id_{\K(L^2(G))})(\cdot)(\id_A\o W)$, with $\sigma : C^*_r(G)\o
\K(L^2(G))\> \K(L^2(G))\o C^*_r(G)$ denoting the swap $\sigma(x\o y)=y\o x$, and with $W$
denoting the unitary $Wf(s,t)=f(s, s^{-1}t)$.

Both for actions and for coactions, under this duality isomorphisms, an equivariant
$*$-homomorphism $\f: A\>B$ becomes the equivariant $*$-homomorphism $\f\o
\id_{\K(L^2(G))}: A\o \K(L^2(G))\> B\o \K(L^2(G))$.

Let $G$ be a locally compact group, and let $A$ and $B$ be $G$-algebras; the
Baaj-Skandalis duality comprises an isomorphism between $KK^G(A,B)$ and
$KK^{C^*_r(G)}(A\r G, B\r G)$, along with a similar result for $\hat G$-algebras, i.e.,
an isomorphism between $KK^{C^*_r(G)}(A,B)$ and $KK^{G}(A\r \hat G, B\r \hat G)$.

We now prove this result in $E$-theory. First we need a lemma.

\begin{lema} For every exact group $G$, and every equivariant short exact sequence
of $\hat G$-algebras $0\>J\stackrel{i}{\>}B\stackrel{p}{\>}A\>0$, the sequence $0\>J\r
\hat G \stackrel{i\times 1}{\>}B\r \hat G\stackrel{p\times 1}{\>}A\r \hat G\>0$ of
reduced cross product by the coactions of $G$, is also exact.
\end{lema}
{\bf Proof.} We have to prove the exactness in the middle, or more precisely the
inclusion $\im(i\times 1)\supseteq \ker(i\times 1)$. Suppose that $\im(i\times 1)$ is a
proper ideal of $\ker(i\times 1)$, and take the cross product by the action of $G$. If
follows  that $\im(i\times 1)\r G$ is a proper ideal of $\ker(p\times 1)\r G$. Identify
$\im(i\times 1)\r G$ with $\im i \o \K(L^2(G))$, and $\ker(p\times 1)\r G$ with $\ker p
\o \K(L^2(G))$, and a contradiction follows. $\square$

\begin{teo} Let $G$ be an exact group, and let $A$ and $B$ be $G$-algebras; there is an
isomorphism
$$E^G(A,B)\simeq E^{C^*_r(G)}(A\r G, B\r G).$$
If $C$ and $D$ are $\hat G$-algebras, there is an isomorphism
$$E^{C^*_r(G)}(C,D)\simeq E^G(C\r \hat G, D\r \hat G).$$
\end{teo}
{\bf Proof.} We shall define  a map $$J_G: \{A,B\}^{\overline{C_0(G)}}\> \{A\r G ,B\r
G\}^{\overline{C_r^*(G)}}$$ and a map  $$J_{\hat G}: \{C,D\}^{\overline{C_r^*(G)}}\>
\{A\r \hat G ,B\r \hat G\}^{\overline{C_0(G)}},$$ such that $J_{\hat G}\circ J_G (\f)=
\f\o \id_{\K(L^2(G))}$ for every asymptotic morphism $\f\in \{A,B\}^{\overline{C_0(G)}}$
and such that $J_{G}\circ J_{\hat G} (\psi)= \psi\o \id_{\K(L^2(G))}$ for every
asymptotic morphism $\psi\in \{C,D\}^{\overline{C_r^*(G)}}$.

First, observe that if $\alpha$ and $\beta$ are two exterior equivalent actions of $G$ on
an algebra $A$, then the corresponding cross product $A\r G$ and $A\r G$ are isomorphic.
A similar property holds for coactions of groups. This can be proved either by a direct
calculation, or be seen as a particular case of \cite[Proposition 7.6]{BS|unitmult}.

If follows that cross product is a functor $\cdot \r G:
\overline{C_0(G)}\mathrm{-alg}\>\overline{C_r^*(G)}\mathrm{-alg}.$ Also, the $\hat
G$-algebras give rise to a subcategory of $\overline{C_r^*(G)}\mathrm{-alg}$, a category
that we denote $\overline{\hat G}\mathrm{-alg}$, and the cross product by coactions of
$G$, gives a functor $\cdot \r \hat G: \overline{\hat
G}\mathrm{-alg}\>\overline{C_0(G)}\mathrm{-alg}.$

The question of which sufficient condition a functor $F$ has to fulfill in order to have
a natural extension to the category of homotopy classes of asymptotic morphisms, has been
studied in \cite[Chapter 3]{GHT}. Exactness plays an important role, and the argument we
used in Lemma 2.5 is quite similar to the general setting quoted above.

Thanks to the exactness of $G$, and the lemma above, if follows the existence of $J_G$
and $J_{\hat G}$ extending the functors cross product.

The condition $J_{\hat G}\circ J_{G} (\f)= \f\o \id_{\K(L^2(G))}$ for every $\f\in
\{A,B\}^{\overline{C_0(G)}}$, holds if $\f$ is a $*$-homomorphism. Moreover, we can
regard this equality as a natural transformation of two functors from the category
$\overline{C_0(G)}\mathrm{-alg}$ to itself. Extension of functors to asymptotic morphisms
comes along with extension of natural transformation between them, as stated by
\cite[Proposition 3.6]{GHT}. Hence this equality holds for every $\f\in
\{A,B\}^{\overline{C_0(G)}}$. A similar argument applies for the composition $J_{G}\circ
J_{\hat G}$. The maps $J_G$ and $J_{\hat G}$ induce maps in $E$-theory $E^G(A,B)\>
E^{C^*_r(G)}(A\r G, B\r G)$ and $E^{C^*_r(G)}(C,D)\simeq E^G(C\r \hat G, D\r \hat G)$
which are inverse to each other. $\square$

\begin{rmq}\em As explained in \cite[Remark 7.7(b)]{BS|unitmult}, the duality isomorphism
holds in a more general setting, that of a pair of dual Hopf $C^*$-algebras associated to
a multiplicative unitary. The proofs of the theorem and of the lemma preceding it adapts
to this case, provided that we assume appropriate exactness condition.
\end{rmq}

\section{The universal property of $KK^S$}

In this section we prove the existence of a natural transformation from equivariant
$KK$-theory to equivariant $E$-theory. To this end, we prove an universal property which
characterizes $KK$-theory as a category which is stable, homotopy invariant and split
exact. A similar result was proved in the non equivariant case by Higson in \cite{Hig2}
based on previous work by Cuntz, and in the group case by Thomsen in
\cite{Thomsen|univKK}. The proof is based on a description of the $KK$-theory in terms of
quasi-morphisms, commonly referred as the Cuntz' picture.

Let us first remind the notion of a Hilbert module endowed with a coaction of a Hopf
$C^*$-algebra. Let $S$ be a Hopf $C^*$-algebra, and assume that $A$ is an $S$-algebra
with coaction $\a : A\> \M(A\o S)$ and let $E$ be a Hilbert $A$-module. Identify the
Hilbert $A\o S$-module $E\o S$ with $\K(A\o S, E\o S)$ and denote by
$$\M(E\o S)= \{T\in \L(A\o S,E\o S) |\ \forall s\in S \ (1_E\o s)T \ \mathrm{ and }\ \
T(1_B\o s) \in E\o S \},$$ the $S$-multipliers of the module $E$ (see section 1.4 from
\cite{EKQR} for more details on this construction). It is a  Hilbert $\M(A\o S)$-module
with scalar product given by $<T_1,T_2>=T_1^*T_2 \in \M(A\o S) \subset \L(A\o S)$.
 A coaction of $S$ on $E_A$ is given by a linear map $\d_E :
E\> \M(E\o S)$ such that
\begin{enumerate}
\item $\d_E(e)\a(a)=\d_E(ea)$,
\item $\a(<e,f>)=<\d_E(e),\d_E(f)>$ for all $a\in A , e,f \in E$,
\item $\d_E(E)(A\o S)$ is dense in $E\o S$ and
\item $(\id_E\o \d_A)\circ \d_A = (\d_E\o \id_S)\circ \d_A$.
\end{enumerate}

The last condition is about the extensions of these maps to $\L(A\o S, E\o S)\>\L(A\o S\o
S,E \o S \o S)$. A second way of defining a coaction is like follows: denote by $T_e\in
\L(A\o S, E\o_{\a} (A\o S))$ the operator given by $T_e(x)=e\o_{\a} x$ for $x\in A\o S$,
and for every $e \in E$. An unitary $V\in \L(E\o_{\a} (A\o S), E\o S)$ is {\it
admissible} if for every $e\in E$, $VT_e\in \M(E\o S)$ and $(V\o \id_S)(V\o_{\a\o \id_S}
1)= V\o_{\id_A\o \d_S}1\in \L(E\o_{\a^2}(A\o S\o S), E\o S\o S)$. There is a one-to-one
correspondence between the coactions on $E$ and the admissible unitaries
\cite[Proposition 2.4]{BS}.

Assume that  $A$ and $B$ are $S$-algebras, and that $E$ is a Hilbert $S$-$B$-module; an
$*$-homomorphism $\pi:A\>\L(E)$ is $S$-equivariant if for all $a\in A$ and $e\in E$,
$\d_E(\pi(a)e)=(\pi\o \id_S)(\a)(a)\circ \d_E(e).$

Let us first recall the definition of the equivariant $KK$-theory of Baaj and Skandalis.
Let $A,B$ be two graded $S$-algebras. A Kasparov triple $(E,\pi,T)$ is given by a Hilbert
$S$-$B$-module $E$, a grading preserving $S$-equivariant representation $\pi:A\>\L(E)$,
and a degree one operator $T\in\L(E)$ such that
\begin{itemize}
\item[(i)]
$\pi(a)(T-T^*)\in\K(E)$;
\item[(ii)]
$\pi(a)(T^2-1)\in\K(E)$;
\item[(iii)]
$[T,\pi(a)]\in\K(E)$;
\item[(iv)]
for all $x \in A\o S$, $(\pi \o \id_S)(x)(F\o_\mathbb{C} 1- V(F\hat\o_{\b} 1)V^*)\in
\K(E\o S)$, where $V\in \L(E\o_{\b} (B\o S), E\o S)$ is the unitary defining the coaction
of $S$ on the Hilbert module $E$.
\end{itemize}

The last condition can also be written as follows $(F\o_\mathbb{C} 1- V(F\hat\o_{\b}
1)V^*) (\pi \o
\id_S)(x)\in \K(E\o S)$, for all $x \in A\o S$.\\
Such a triple is called degenerate if $\pi(a)(T-T^*)=\pi(a)(T^2-1)=[T,\pi(a)]=0$ for all
$a\in A$ and if $(\pi \o \id_S)(x)(F\o_\mathbb{C} 1- V(F\hat\o_{\b} 1)V^*)=0$, for all $x
\in A\o S$.  The set of unitary equivalence classes of Kasparov triples form a semigroup.
The group $KK^S(A,B)$ is defined as the quotient of this semigroup by a homotopy relation
which is defined using triples for the pair of algebras $(A, B[0,1])$. Degenerate triples
are homotopic to zero.

From now on, algebras are trivially graded.

\begin{defi} \em An $S$-quasi-morphism is a pair $(\f_+,\f_-):A\>M(B\otimes\K)$ of
quasi-unital morphisms, along with a pair $(V_+,V_-)\in M(B\o \K \o S)$ of $\b\o
\id_\K$-cocycles such that
\begin{itemize}
\item[(i)]
$\f_+(a)-\f_-(a)\in B\otimes\K,$ for all $a\in A;$
\item[(ii)]
$\f_\pm$ are $\d_A$-$\d_{B\o\K}^\pm$-equivariant where $\d_{B\o\K}^\pm$ denote the
coaction $\d_{B\o\K}^{\pm}(\cdot):=V_\pm(\b(\cdot)\o \id_\K)V^*_\pm$;
\item[(iii)]
$V_+ -V_- \in \M(B\o\K\o S)$.
\end{itemize}
\end{defi}

An $S$-quasi-morphism is called {\it degenerate} if $\f_+=\f_-=0$.

Two $S$-quasi-morphisms $(\f_\pm,U_\pm)$ and $(\psi_\pm,V_\pm)$ are {\it isomorphic} if
there is a $\b\o \id_\K$-equivariant automorphism  $\Theta$ of $B\o \K$ such that
$\psi_\pm =\Theta \circ \f_\pm$ and such that $V_\pm=(\Theta \o \id_S)\circ U_\pm$. Two
$S$-quasi-morphisms $(\f_\pm^0,U_\pm^0)$ and $(\f_\pm^1,U_\pm^1)$ are {\it homotopic} if
there is an $S$-quasi-morphism $(\Phi_\pm,U_\pm):A\>M(B\o \K \o C[0,1])$ such that
$$(\f_\pm^i,U_\pm^i)= (ev_i\circ\Phi_\pm,(ev_i\o\id_S)\circ U_\pm),$$
where $ev_i:M(B\o\K\o C[0,1])\> M(B\o \K)$ denote the evaluations in $i=0$ and in $i=1$.
The sum the $S$-quasi-morphisms $(\f_\pm,U_\pm)$ and $(\psi_\pm,V_\pm)$ is defined using
an isomorphism $M_2(\K)\simeq \K$ as $(\f_\pm\oplus\psi_\pm ,U_\pm\oplus V_\pm)$; its
homotopy class does not depend on the choice made. Finally, two $S$-quasi-morphisms
$(\f_\pm,U_\pm)$ and $(\psi_\pm,V_\pm)$ are {\it equivalent} if there are degenerate
$S$-quasi-morphisms $(\f_\pm ' ,U_\pm ')$ and $(\psi_\pm ',V_\pm ')$ such that
$(\f_\pm,U_\pm)\oplus (\f_\pm ' ,U_\pm ') $ and $(\psi_\pm,V_\pm)\oplus (\psi_\pm ',V_\pm
')$ are homotopic.

\begin{teo}
$KK^S(A,B)$ is isomorphic with the group of equivalence classes of $S$-quasi-morphisms
from $A$ to $B$.
\end{teo}
{\bf Proof.} We assume that $B$ is trivially stable, this in order to eliminate from what
follows an extra copy of $\K$, with the trivial coaction on it, this when we think that
this will not hinder the argument. Let $(\f_\pm,V_\pm)$ be a $S$-quasi-morphism, we
associate to it the following triple
$(E,\pi ,T)$:\\
- the Hilbert module $E$ is the direct sum $B\o\K \oplus B\o\K$ with the coaction of
$S$ given by $V_+(\b\o \id_{\K})\oplus V_-(\b\o \id_{\K})$. \\
- the representation $\pi:A\>\L_{B\otimes\K}(E)$ and the operator $F\in
\L_{B\otimes\K}(E)$ respectively defined by
$$\pi(a)=\left(\begin{array}{cc}
\f_+(a) & 0\\
0 & \f_-(a)\end{array}\right)\ \mbox{and by}\ F=\left(\begin{array}{cc}
0 & 1\\
1 & 0\end{array}\right).$$

Note that the condition $V_+ -V_- \in \M(B\o\K\o S)$ writes as the equivariance condition
for this particular triple. In this way one associates a degenerate triple to a
degenerate $S$-quasi-morphism; it agrees with the direct sums and with the homotopies of
$S$-quasi-morphisms and of equivariant Kasparov's triples. Hence it defines a map from
the set of equivalence classes of $S$-quasi-morphisms $(\f_\pm,V_\pm):A\>M(B\otimes\K)$
to the group $KK^S(A,B)$.

Conversely, take $(E,\pi,T)\in KK^S(A,B\o\K)$. One can assume that the representation
$\pi$ is nondegenerate, as the argument from \cite[Proposition 18.3.6]{B} applies to the
equivariant case too. The idea of the proof is first to replace the Hilbert module $E$
with the Hilbert module $B\o\K\oplus B\o\K$ and then to replace the operator $T$ with the
operator  $\left(\begin{array}{cc}
0 & 1\\
1 & 0\end{array}\right)$. In then nonequivariant case, this is described in detail
\cite[Section 17.6]{B} and in the group case it is done by Thomsen in
\cite{Thomsen|univKK}. We describe how this transformations work in our case.

 Denote by $V:E\>\M(E\o S)$ the unitary implementing the coaction of
$S$ on the Hilbert module $E$. The coaction on the $B\o\K$-module $B\o\K\oplus B\o\K$ is
the standard one.

The triple $(B\o\K \oplus B\o\K,0,0)\in KK^S(A,B\otimes\K)$ is degenerate thus
$$(E,\pi,T)\oplus(B\o\K \oplus B\o\K,0,0)=(E,\pi,T).$$

The Kasparov stabilization theorem states that there is a non equivariant graded
isomorphism of Hilbert $B\o \K$-modules $$ E\oplus l^2(B\o \K)\oplus l^2(B\o \K)\simeq
l^2(B\o \K)\oplus l^2(B\o \K).$$ Moreover, because $B\o \K$ is stable, there is a graded
isomorphism of Hilbert $B\o \K$-modules between $l^2(B\o \K)\oplus l^2(B\o \K)$ and $B\o
\K\oplus B\o \K$ as shown in \cite[Lemma 1.3.2]{JT}; denote by
$$\Psi: E\oplus  B\o \K\oplus B\o \K\> B\o \K\oplus B\o \K$$
the resulting isomorphism. Let $W$ be the $*$-homomorphism making the following diagram
commutative:
$$
\begin{diagram}
\node{E\oplus (B\o \K \oplus B\o \K)}\arrow{e,t}{\Psi} \arrow{s,l}{V \oplus \b\o
\id_\K\oplus\b\o \id_\K}
\node{B\o \K\oplus  B\o \K}\arrow{s,r}{W}\\
\node{\M(E\o S\oplus ( B\o \K\oplus  B\o \K)\o S)}\arrow{e,t}{\Psi\o \id_S} \node{\M((B\o
\K\oplus B\o \K)\o S);}
\end{diagram}
$$
it is a coaction on $B\o \K\oplus B\o \K$, in general different from the standard one, as
we cannot assume that the stabilization morphism $\Psi$ is $S$-equivariant.

This takes care of the Hilbert module. We can further simplify the resulting triple  and
obtain a triple
$$(B\otimes\K\oplus B\otimes\K,\f, \left(\begin{array}{cc}
0 & 1\\
1 & 0\end{array}\right)).$$ Note that the representation $\f$ is quasi-unital.

The action on this module has the form  $V =(V_+,V_-)$ with $V_+,V_-$ denoting the
coactions of $S$ on the Hilbert $B\o\K$-module $B\o\K$. Denote $V_+,V_-\in\L(B\o \K \o
S,B\o\K \o S)$ the unitaries defining this action, and set $U_\pm \in U(M(B\o \K \o S))$,
the $\b\o \id_\K$-cocycles
$$U_\pm =V_\pm \circ(\b\o \id_\K)^*.$$

Take $\f=\left(\begin{array}{cc}
\f_+ & 0\\
0 & \f_-\end{array}\right)$; $\f_\pm$ is $\d_A$-$U_\pm(\d_B\o
\id_\K)U_{\pm}^*$-equivariant. It follows from the condition $[\f(a),F]\in B\otimes\K$
that $\f_+(a)-\f_-(a)\in B\otimes\K$, and it follows from the equivariance condition that
$U_+-U_-\in \M(B\o \K \o S)$, thus $(\f_\pm,U_\pm)$ is an $S$-quasi-morphism.

The rest of the proof, follows like in the group case.
$\square$\\

That $KK^S(A,\cdot)$ is itself a functor which is stable, homotopy invariant, and
split-exact. Stability is best express through Morita equivalence: an equivariant
imprimitivity bimodule provides an invertible element in equivariant $KK$-theory.
Homotopy invariance follows from the definitions. The argument from \cite{Thomsen|univKK}
adapts without changes to prove split-exactness.

\begin{prop}(Universal property of $KK^S$) Let  $F:S\mathrm{-alg}\>Ab$ be a covariant functor
which is stable, homotopy invariant, and split-exact. Then for every $S$-algebra $A$ and
every element $d\in F(A)$ there exists an unique natural transformation
$T_A:KK^S(A,\cdot)\>F(\cdot)$ such that $T_A(1_A)=d$.
\end{prop}
{\bf Proof.} We indicate how an element $x\in KK^S(A,B)$ defines an application
$L(x):F(A)\>F(B)$, as we do it a bit differently from \cite{Thomsen|univKK}. The  rest of
the proof, follows like in the non-equivariant case.

Let $(\f_\pm,V_\pm)$ be a $S$-quasi-morphism representing $x$.

Consider the coaction $\bar \d_B$ on $M_2(B\o \K)$ which on the upper left corner
restricts to $\d_B^{+}$ and on the lower right corner restricts to $\d_B^{-}$. This
action is exterior equivalent to $\d_B$, with cocycle diag$(V_{+},V_{-})$, hence there is
an isomorphism
$$ I_B: F(B,\d_B)\>F(M_2(B\o \K),\bar \d_B).$$
The maps $i_B^{+}(b)= \ma b 0 0 0 $ and  $i_B^{-}(b)= \ma 0 0 0 b $, induce respectively
the isomorphisms $i_{B*}^{+} : F(B\o \K, \d_B^+)\> F(M_2(B\o \K),\bar \d_B)$ and
$i_{B*}^{-} : F(B\o \K, \d_B^-)\> F(M_2(B\o \K),\bar \d_B)$.

Denote by $A_x$ the subalgebra of $A\oplus M_2(M(B\otimes\K))$ of elements $(a,m)$ such
that:\\ $\f_+(a)- m_{11} \in B\otimes\K$, $\f_-(a)- m_{22} \in B\otimes\K$ and  $m_{12},
m_{21} \in B\o \K $.

It is an $S$-algebra with coaction $(\d_A, \bar \d_B)$, and it fits into an equivariant
short exact sequence
$$ 0\longrightarrow (M_2(B\o\K),\bar \d_B)\stackrel{j}{\longrightarrow}A_x
\stackrel{p}{\longrightarrow}(A,\d_A)\longrightarrow 0$$ with maps given by $j(b)=(0,b)$
and by $p(a,m)=a$.

This short exact sequence splits equivariantly into two different ways by the
$*$-homomorphisms $s_\pm:A\>A_x$, given by
$$s_+(a)=\left(a,\ma {\f_+(a)} 0 0 0 \right) \ \mbox {and} \ s_-(a)=\left(a,
\ma 0 0 0 {\f_-(a)}\right). $$

Define $L(x):F(A)\>F(B)$ as the composition
$$ \Phi: F(A)\>F(A_x)\>F(M_2(B\o\K), \bar \d_B)\>F(B,\d_B)$$
given by $\Phi= I_{B *}^{-1}\circ j^{-1}_*\circ (s_{+*}-s_{-*})$. $\square$

\begin{rmq}{\rm Let us now describe the product in $KK^S$ using Cuntz' picture, a question
implicit in \cite{Thomsen|univKK}. Let $A$, $B$, $C$ be $S$-algebras and consider
$(\f_\pm,U_\pm)\in KK^S(A,B)$ and $(\psi_\pm,V_\pm)\in KK^S(A,B)$, the product of the two
elements is given by the $S$-quasimorphism
$$(\phi_+\circ \f_+\oplus \phi_-\circ \f_-, \phi_-\circ \f_+\oplus \phi_+\circ \f_-)$$
with the associated cocycle
$$  (V_+\f_+^\sharp(U_+)\oplus V_-\f_-^\sharp(U_-),V_-\f_+^\sharp(U_+)\oplus
V_+\f_-^\sharp(U_-)).$$ This defines a product in the category $KK^S$, and the claim
follows from the fact that such a product is unique. This in itself is a consequence of
the universal property.}
\end{rmq}

\begin{coro} There is a natural transformation from the category $KK^S$ to the category
$E^S$.
\end{coro}

\section{Appendix: equivariant stabilization}

In this appendix we examine various ways in which stability of a functor can be
expressed.

Let $F:S\mathrm{-alg}\> Ab$ be a  functor from the category of $S$-algebras to abelian
groups. The following definition is the notion of stability used by Thomsen to
characterize $KK$-theory.

\begin{defi} \em The functor $F:S\mathrm{-alg}\> Ab$ is stable if for every $S$-algebra
$(A, \d_A)$ and every compatible action $\d_{A\o\K}$ on $A\o\K$, the compatibility
morphism $a\>a\o e$ induces an isomorphism $F(A,\delta_A)\simeq F(A\o\K,\d_{A\o\K})$.
\end{defi}

Remark that if $F$ is stable and if $V$ is a $\d_A$-cocycle then $F(A,\delta_A)\simeq
F(A,\d_A^V)$, as by the lemma 3.4 they are both compatible with $(M_2(A)\o \K,
\d_{M_2(A)}\o \id_\K)$.\\

In the non equivariant case, stability can be stated in terms of Morita equivalence. Let
us recall this notion: two algebras $A$ and $B$ are Morita equivalent if there is a
Hilbert $B$-module which is full, i.e. such that $<E,E>=B$, and such that $A\simeq
\K_B(E)$; the module $E$ is called imprimitivity bimodule.

The Brown-Green-Rieffel theorem (\cite{MP}) states that two $C^*$-algebras $A$ and $B$
are Morita equivalent if and only if they are stably isomorphic, that is, if $A\o
\K\simeq B\o \K$. Hence, in the non-equivariant case, a functor is stable if and only if
it is Morita invariant in the obvious sense.

The notion of Morita equivalence extends naturally to the equivariant case: for a
Hopf-$C^*$-algebra $S$ and two $S$-algebras  $A$ and $B$ assume that there is a Hilbert
$S$-$B$-module $E$ which is full and such that the isomorphism $A\simeq \K_B(E)$ is $S$
-equivariant.\\
\\
{\bf Examples.}1. Let $G$ be a group and $(A,\alpha)$ be a $G$-algebra; then $(A,\alpha)$
and $(A\o \K(L^2(G)),\alpha\o \lambda_G)$ are $G$-Morita equivalent with imprimitivity
bimodule $L^2(G,A)$.
\\
2. Let $(A, \delta_A)$ be an $S$-algebra; it is Morita equivalent to $(A, \delta_A')$ if
and only if the actions  $\delta_A$ and $\delta_A'$ are exterior equivalent. To see this,
note that the linking algebra is in this case $M_2(A)$ and that the action on diagonal
elements is the diagonal $\small \ma {\d_A} 0 0 {\d'_A}$, hence the claim follows from
the lemma 3.4.

The imprimitivity bimodule is the algebra $A$ itself, seen as a Hilbert $A$-module,
endowed with the coaction $\d: A\>\M(A\o S)$, $\d(a)= V\d_A(a)$, where $V$ denotes the
$\d_A$-cocycle for which $\d_A'(\cdot)= V\d_A(\cdot)V^*$.
\\
3. Let $(A, \delta_A)$ be an $S$-algebra; a coaction $(A\o\K, \d_{A\o\K})$ is compatible
with $\d_A$ if and only if the two algebras are Morita equivalent through an
imprimitivity bimodule $A\oplus A^\infty_A$ with the coaction on the first factor $A$
given by $\d_A$. The nontrivial implication follows from \cite[Proposition 2.7 (a)]{BS}.
\\
\\
Even though not needed here, it is worth mentioning that Morita equivalence appears as
more natural from the point of view of the categories involved (see \cite{EKQR}).

The following is an equivariant version of the Brown-Green-Rieffel theorem, which shows
that the notion of stability we use, is also equivalent to equivariant Morita
equivalence. In the group case, it is due to Combes \cite{combes|crossed products}; his
proof is different though.

\begin{teo} Let $(A,\delta_A)$ and $(B,\delta_B)$ be $S$-algebras, $S$-equivarianly
Morita equivalent. There is a coaction $\delta_{B\o\K}$ of $S$ on $B\o\K$, which is
compatible with $\delta_B$, such that the algebras $(A\o\K, \delta_A\o \id_\K)$ and $(B\o
\K,\delta_{B\o\K})$ are $S$-isomorphic.
\end{teo}
{\bf Proof.} Let us first remaind that the proof of the (non-equivariant)
Brown-Green-Rieffel theorem is based on the following result: for every full Hilbert
$B$-module, there is an isomorphism of Hilbert $B$-modules between $E^\infty$ and
$B^\infty$. If follows  that if $A$ and $B$ are Morita equivalent with imprimitivity
bimodule $E$, then
$$A\o \K\simeq \K_B(E^\infty)\simeq \K_B(B^\infty)\simeq B\o\K. $$

Assume that $A$, $B$ and $E$ are endowed with coactions of $S$, denoted respectively by
$\delta_A$, $\delta_B$, and $\delta_E$. The first isomorphism above is equivariant for
the coaction $\delta_A\o \id_\K$ on $A\o\K$ and the coaction on $\K_B(E^\infty)\simeq
\K_B(E)\o \K$ induced by the coaction $\delta_E$ on $E$ and the trivial coaction on $\K$.

Consider the map $\beta$ making the following diagram commutative:
$$
\begin{diagram}
\node{E^\infty}\arrow{e,t}{T}\arrow{s,l}{\delta_E^\infty} \node{B^\infty}
\arrow{s,r}{\beta}\\
\node{\M(E^\infty\o S)}\arrow{e,t}{T\o\id_S}\node{\M(B^\infty\o S)}
\end{diagram}
$$
where $T:E^\infty\>B^\infty$ denotes a (non-equivariant) isomorphism of Hilbert
$B$-modules. It is easy to check that $\beta$ defines a coaction on the Hilbert
$S$-$B$-module $B^\infty$, and we claim that it induces on the algebra
$\K_B(B^\infty)\simeq B\o\K$ a coaction $\delta_{B\o\K}$ which is compatible with
$\delta_B$, and for which the above isomorphism $A\o\K\simeq B\o\K$ is equivariant.
Consider the isomorphism of Hilbert $B$-modules
$$T\oplus\id_{B^\infty}:E^\infty\oplus B^\infty\>B^\infty\oplus B^\infty;$$
it is $\d_E^\infty\oplus\d_B^\infty$-$\beta\oplus\d_B^\infty$-equivariant. This
isomorphism induces an isomorphism of $S$-algebras between the corresponding algebras of
compact operators, and hence there is an action on $\K(E^\infty\oplus B^\infty)\simeq
M_2(B\o\K)$ which on left-upper corner is $\d_{B\o \K}$ and on the lower-right corner is
$\delta_B\o \id_\K$, which proves compatibility. The rest is obvious. $\square$

To summarize the discussion above, stabilization can be define in three
equivalent ways:\\
1. (stable) $F(A,\delta_A)\simeq F(A\o\K,\d_{A\o\K})$ for every $S$-algebra
$(A,\delta_A)$ and every compatible coaction $\d_{A\o\K}$ on $A\o\K$,\\
2. (Morita equivalence) $F(A,\delta_A)\simeq F(B,\d_B)$ for every $S$-Morita equivalent
pair of $S$-algebras $(A,\delta_A)$ and $(B,\d_B)$, and\\
3. (exterior equivalence along with trivial stability) $F(A,\delta_A)\simeq F(A,\d_A')$
for a pair of exterior equivalent coactions, and $F(A,\delta_A)\simeq
F(A\o\K,\d_A\o\id_{K})$.

Let us now take a closer look at the case of a group $G$ action.

\begin{lema} A functor $F$ defined on the category of $G$-algebras is stable if and only
if $F(\cdot )\simeq F(\cdot \o \K_G)$.
\end{lema}
{\bf Proof.} Let $A$ be a $G$-algebra; $A$ and $A\o \K_G$ are $G$-Morita equivalent so
the condition is necessary for stability.

Assume now that $F(A,\alpha)\simeq F(A\o \K_G, \alpha \o \lambda_G)$ for every action
$\alpha$ of $G$ on $A$. Note first that $(\K_G, \lambda_G)\simeq (\K_G\o \K , \lambda_G\o
\id_{\K})$, hence we only need to check that if $\alpha$ and $\alpha'$ are two exterior
actions of $G$ on $A$, then $F(A,\alpha)\simeq F(A,\alpha')$. Write $\alpha'_g(\cdot)=
u_g \alpha_g(\cdot)u^*_g$; and consider the $(A,\alpha')$-$(A,\alpha)$-imprimitivity
bimodule $A$, with its coaction $\gamma$ given by $\gamma_g(a)=u_g\alpha_g(a)$, for $a\in
A $ and $g\in G$. The Hilbert $G$-$A$-modules $(A_A,\gamma)$ and $(A_A, \alpha)$ are non
equivariantly isomorphic, thus $(L^2(G,A)_A,\lambda_G \o \gamma)$ and
$(L^2(G,A)_A,\lambda_G\o \alpha)$ are $G$-isomorphic. Take their compact operators, which
are respectively $(A\o \K_G, \alpha'\o \lambda_G)$ and $(A\o \K_G, \alpha\o \lambda_G)$,
are $G$-isomorphic, hence $F(A,\alpha)\simeq F(A,\alpha')$.$\square$

If follows that given a functor $F$, one can replace it with the functor $\bar F(\cdot):=
F(\cdot \o \K_G)$, which is stable.

{\sc Departamento de Matem\'atica, Instituto Superior T\'ecnico, Av. Rovisco Pais,
1049-001, LISBOA, Portugal, popescu@math.ist.utl.pt}
\end{document}